\documentclass[modern]{aastex631_cranmer}
\usepackage{hyperref}
\usepackage{epsfig}
\usepackage{amsmath,amsfonts,MnSymbol,amssymb,mathrsfs}
\usepackage{cjhebrew}
\newcommand{\expon}{\,\, \mbox{\textasteriskcentered \!\! \textasteriskcentered \,}}

\shorttitle{Intuitive but Non-Rigorous Explanations of Infinite Numbers}
\shortauthors{S.\  R.\  Cranmer}

\begin{document}

\title{Intuitive Explanations of Infinite Numbers for Non-Specialists}

\correspondingauthor{Steven R. Cranmer}
\author[0000-0002-3699-3134]{Steven R. Cranmer}
\affiliation{Department of Astrophysical and Planetary Sciences,
Laboratory for Atmospheric and Space Physics,
University of Colorado, Boulder, CO, USA}

\begin{abstract}
The mathematical study of infinity seems to have the ability to
transport the mind to lofty and unusual realms.
Decades ago, I was transported in this way by Rudy Rucker's book
{\em Infinity and the Mind.}
Despite much subsequent learning and teaching of mathematics
in the service of physics and astronomy, there remained quite a few
aspects of the ``higher infinities'' that I was still far from
comprehending.
Thus, I wanted to dive back in to understand those ideas and
to find good ways to explain them to other non-specialists.
These notes are an attempt to do this.
They begin by discussing huge (but still describable) finite
numbers, then they proceed to the concepts of fixed points,
Cantor's countably infinite ordinals, transfinite cardinals and
the continuum hypothesis, and then to more recent attempts to
define still-larger infinities called large cardinals.
I should warn the reader that the author is an astrophysicist, so
these pedagogic explanations may not be satisfying (or anywhere
nearly sufficiently rigorous) to a mathematician.
Still, I hope these explanations provide some intuitive insights
about concepts that are too large to fit in the physical universe.
\end{abstract}

\section{Introduction}
\label{sec:intro}

These notes are an attempt to describe an allegedly ``useless''
field of pure mathematics.
My own experience, via research and teaching, has been much more
pragmatic, with a strong emphasis on the ``unreasonable effectiveness
of mathematics'' \citep{Wi60} as applied to the physical world.
Why, then, focus on the mathematical study of infinity?
Its odd reputation precedes us by more than a century.
In response to Georg Cantor's development of infinite set theory in
the 19th century, fellow mathematician Leopold Kronecker called those ideas
``mathematical nonsense,'' and in turn called Cantor a ``scientific
charlatan'' and a ``corruptor of youth.''
Henri Poincar\'{e} called these theories ``a grave disease infecting
the body of mathematics.''
Extreme skepticism about the validity and usefulness of infinite
mathematics continues to the present \citep[e.g.,][]{Zn04,Muck23}.

On the other hand, David Hilbert lauded this field of study as
``one of the supreme achievements of purely intellectual human activity,''
and also remained steadfast by saying that ``no one shall be able to expel
us from the paradise that Cantor has created for us'' \citep{Hi26}.
Building on the mystical poetry of William Blake, \citet{Jg13}
asserted that ``It is through mathematics that we can hold not one but an
infinity of infinities in the palm of our hand, for all of eternity.''
There is beauty and wonder here, and these notes are an attempt to
convey that beauty and wonder to others.\footnote{%
I have no idea whether the process of reading about and understanding
these concepts could assist the reader in approaching a mystical
or transcendent state, as \citet{Ru82} suggested, but I hope that may
be possible.  Thus, an alternate title for these notes could be
``Grokking the Infinite.''}

It should be noted that the process of developing rules and theorems
about infinity does seem to have been useful to obtaining results in
other fields of mathematics \citep[e.g.,][]{Dh95,Pe10}
and computing \citep{Da18},
and it has helped improve our understanding of some aspects of physics
\citep{Aug84,NS99,TW01,EN02}.
However, this usefulness sidesteps the question about whether or
not the actual infinities of mathematics ``actually'' exist in our own
physical universe.
Many books on infinity often spend some time talking about cosmology,
which we may define as the study of the {\em utter vastness} of our
universe in both space and time.
Despite my admiration for Giordano Bruno,
I really have no idea if there are any modern theories of cosmology
that require the universe to have a truly infinite extent in
either space or time \citep[see, perhaps,][]{HW11}.
These notes are not about that.

I should probably mention a few other things that these notes
do not cover.
I will assume that the standard use of infinity (as a limit)
in calculus and mathematical physics is reasonably well-covered by
others.
Thus, we won't need to explore any deep mysteriousness in expressions
such as
\begin{equation}
  \int_{-\infty}^{+\infty} dx \,\, e^{-x^2} \, = \, \sqrt{\pi}
  \,\,\,\,\,\,\,\,\,\,\,\,\,\,\,\,
  \mbox{or}
  \,\,\,\,\,\,\,\,\,\,\,\,\,\,\,\,
  \lim_{n \rightarrow \infty} \left( 1 + \frac{x}{n} \right)^n
  \,\, = \,\, e^x \,\,\, ,
\end{equation}
despite their immense usefulness.
Also, before anyone brings it up, I think I will refuse to wade
into the clickbait that is
\begin{equation}
  1 \,\, + \,\,
  2 \,\, + \,\,
  3 \,\, + \,\,
  4 \,\, + \,\, \cdots \,\,\, = \,\,\,
  \zeta(-1) \,\,\, = \,\,\, -\tfrac{1}{12} \,\,\, .
\end{equation}
If the above is new and meaningless to you, consider yourself lucky.
However, if you do want to learn more about this bit of mathematical
skullduggery, see \citet{Hr49} or
\href{https://en.wikipedia.org/wiki/1_%2B_2_%2B_3_%2B_4_%2B_%E2%8B%AF}%
{wikipedia}.

What do we actually cover?
Section~\ref{sec:finite} gets us started by contemplating huge,
but still finite numbers.
Section~\ref{sec:powertower} dips a toe into the infinite waters
by exploring what happens to some finite numbers when they are
exponentiated an infinite number of times.
Section~\ref{sec:countable} then makes the leap to the first type of
``countable infinity'' and shows how we can keep counting past
that point.
Section~\ref{sec:real} introduces the continuum of real numbers,
and eases us into the idea that the size of the set of all real numbers
is really and truly {\em bigger} than the set of all integers.
Section~\ref{sec:uncountable} then wades deeper into the pool of
uncountable infinities, and
Section~\ref{sec:beyondC} pauses to take stock of where we stand
after introducing such strange concepts.
Section~\ref{sec:large} then barrels forth to introduce even larger
infinities that are far more ``inaccessible'' than the infinities
that came before.
Section~\ref{sec:morelarge} continues climbing the ladder and
explains how these levels can be understood on the basis of how
complex are the languages used to describe them.
Section~\ref{sec:evenhigher} makes a few final stops before we
run out of oxygen, and
Section~\ref{sec:conc} concludes with some philosophical musings.

\section{Very Big Finite Numbers}
\label{sec:finite}

I can start with my own first exposure to the ``gates'' of infinity:
huge, but still finite, numbers.
I don't remember how old I was when I learned that the pattern in the
names of million ($10^6$), billion ($10^9$), trillion ($10^{12}$),
quadrillion ($10^{15}$), and so on,
could be generalized with Latin prefixes as
\begin{displaymath}
  \mbox{$n$-illion} \,\, \longleftrightarrow \,\,
  10^{3n+3}  \,\,\, .
\end{displaymath}
In the 1970s and 1980s,
annual editions of the {\em Guinness Book of World Records}
had a section on huge numbers that revealed that the $n$-illion
names could be extended at least up to $n=100$ (centillion).
\citet{CG96} described a system that extends them up to $n=1000$
and beyond.
{\em Guinness} also introduced me to
Edward Kasner's googol ($10^{100}$, or 10 duo-trigintillion) and
googolplex ($10^{\rm googol}$), the
Buddhist {\em asankhyeya} ($10^{140}$, or 100 quinto-quadragintillion,
in some translations), as well as more esoteric quantities such as
Skewes's number and Graham's number (see below).

Later, I learned about other ways to reach huge numbers by means of
repeated mathematical operations.
We know that multiplication is repeated addition,
\begin{equation}
  x \cdot y \,\, = \,\,
  \underbrace{x \, + \, x \, + \, x \, + \,
  \ldots \, + \, x}_{y \,\, {\rm times}}
  \label{eq:multdef}
\end{equation}
and that exponentiation is repeated multiplication,
\begin{equation}   
  x^y \,\, = \,\,
  \underbrace{x \, \cdot \, x \, \cdot \, x \, \cdot \,
  \ldots \, \cdot \, x}_{y \,\, {\rm times}}
  \,\,\, .
  \label{eq:expdef}
\end{equation}   
However, what happens when we extend this concept?
The next step would be repeated exponentiation,
\begin{equation}   
  {^y}x \,\, = \,\,
  x^{x^{x^{\udots^{x}}}} \,\,\,\,\,\,\,\, \mbox{($y$ times)}
\end{equation}   
or, in a more computer-friendly notation,
\begin{equation}   
  {^y}x \,\, = \,\,
  \underbrace{x \expon  ( x \expon
  ( x \expon  ( \,\, \cdots \,\,\, (
  x \expon  (x \expon  x) ) ) ) )}_{y \,\, {\rm times}} \,\,\, ,
\end{equation}
where the double-asterisk notation for exponentiation is used by
both Python and Fortran, and the ordering of the parentheses is
important.

Although this type of repeated exponentiation was first discussed
by the Marquis de Condorcet in the 1770s,
it was \citet{Mau01} who devised the ${^y}x$ notation given above,
and \cite{Good47} who called this operation by its (now most
well-known) name of {\em tetration,} the fourth operation in the
series starting with addition.
On a personal note, learning about tetration led me to Lambert's $W$
function \citep{Cor96}, which allowed me to solve some interesting
problems in solar physics \citep[see, e.g.,][]{Cr04}.

\citet{Kn76} created an ``up-arrow notation'' by expressing
exponentiation as $x \uparrow y$ and
tetration as $x \uparrow \uparrow y$.
Thus, taking a number $x$ and tetrating it with itself $y$ times
would be called ``pentation'' ($x \uparrow \uparrow \uparrow y$).
Repeatedly iterated pentation would be ``hexation''
($x \uparrow \uparrow \uparrow \uparrow y$),
repeatedly iterated hexation would be ``septation''
($x \uparrow \uparrow \uparrow \uparrow \uparrow y$), and so on.

There are many ways to further generalize the up-arrow notation
(i.e., Ackermann functions, Conway's chained-arrow notation,
Goodstein's hyperoperator notation, and so on), but the simplest
appears to be just $a[n]b$, where $n$ denotes the order-number of
the operation type: $n=1$ is addition ($a+b$),
$n=2$ is multiplication ($a \cdot b$),
$n=3$ is exponentiation ($a \uparrow b$), $n=4$ is tetration, and so on.
Note that for $n \geq 3$, $a[n]b$ is equal to $a$, followed by $n-2$
up-arrows, followed by $b$.
To illustrate how fast things grow when increasing $n$, we can give
an example:
\begin{eqnarray}
  4 + 4 \, = \, 4[1]4 &\,\,\, = \,\,\,& 8 \nonumber \\
  4 \cdot 4 \, = \, 4[2]4 &\,\,\, = \,\,\,& 16 \nonumber \\
  4 \uparrow 4 \, = \, 4[3]4 &\,\,\, = \,\,\,& 256 \nonumber \\
  4 \uparrow \uparrow 4 \, = \, 4[4]4 &\,\,\, = \,\,\,&
    \mbox{a number with about $10^{154}$ digits} \nonumber  \\
  4 \uparrow \uparrow \uparrow 4 \, = \, 4[5]4 &\,\,\, = \,\,\,&
    \mbox{a number too large to be written with any usual notation.}
    \nonumber
\end{eqnarray}

To explore ever-bigger numbers, we need to shift from describing
specific mathematical
operations to {\em algorithms} that describe iterated operations.
A famous example is Graham's number, which was devised to prove a
theorem in combinatoric number theory \citep[see, e.g.,][]{Gard77}.
Graham's number is defined as $g_{64}$, using a notation defined as
\begin{equation}
  g_n \,\, = \,\, \left\{
  \begin{array}{ll}
    3[6]3 \,\, , & \mbox{for $n=1$,} \\
    3[g_{n-1}+2]3 \,\, , & \mbox{for $n \geq 2$ .} 
  \end{array} \right.
\end{equation}
Thus, $g_1$ (which is defined as 3~$\uparrow\uparrow\uparrow\uparrow$~3)
is essentially indescribable using ordinary numbers, and
$g_2$ (which is writable as 3, followed by
3~$\uparrow\uparrow\uparrow\uparrow$~3 up-arrows, followed by 3)
is even more so.
Going all the way to $g_{64}$ seems somewhat inconceivable, but
the algorithm is straightforwardly defined.

Feel free to take a breath or two at this point.
These concepts can be somewhat dizzying.
However, we {\em can} go further.
Note that our definition of the algorithm that generates $g_{64}$
was short enough to fit into a single paragraph.

\citet{Rayo19} described a competition in which the participants
tried to outdo one another in finding compact ways of describing
the largest possible numbers.
The winner, which is now called {\em Rayo's number,} is essentially
just a description of the {\em size} of the algorithm needed to
describe the number; i.e.,
\begin{quotation}
\noindent
  Rayo's number is the smallest number that is larger than any
  finite number that can be named in the language of first-order
  set theory, using a googol symbols or less.
\end{quotation}

Putting aside the exact meaning of a first-order theory
\citep[see, e.g.,][]{Smith22}---and what one is allowed to say
using the symbols of that language---we would like to know
how we might continue making even bigger numbers.
Let us define the Rayo number as part of an iterated sequence.
We can then call a googol the ``zeroth'' Rayo number
(i.e., ${\cal R}(0) = 10^{100}$).
Thus, if we make the following definition below, Rayo's number
itself is ${\cal R}(1)$.
The definition is:
\begin{quotation}
\noindent
  For $n \geq 1$, we define ${\cal R}(n)$ as
  the smallest number that is larger than any
  finite number that can be named in the language of first-order
  set theory, using ${\cal R}(n-1)$ symbols or less.
\end{quotation}
We now have an iterative algorithm that describes iterative algorithms.
${\cal R}(2)$ must exceed any finite number describable in a number
of symbols given by Rayo's original number,
and it is difficult to fathom how much bigger ${\cal R}(2)$ 
is, in comparison to ${\cal R}(1)$.
What, then, about ${\cal R}(3)$, or ${\cal R}(100)$, or even
${\cal R}(\mbox{googol})$?

Then, of course, we could discuss ${\cal R}({\cal R}(1))$.
This is kind of a meta-Rayo's number.
Once we have realized that doing this is possible, we can define
${\cal R}({\cal R}({\cal R}(1)))$, and then
${\cal R}({\cal R}({\cal R}({\cal R}(1))))$, and so on.
If this is the road we're going down, it would be convenient to define
a new function ${\cal Q}(n,m)$, where $n$ is the number of times that
${\cal R}()$ operates on an original instance of ${\cal R}(m)$.
The three examples at the start of this paragraph would then be
writable more compactly as ${\cal Q}(1,1)$,
${\cal Q}(2,1)$, and ${\cal Q}(3,1)$, respectively.

Well then, what's stopping us from considering numbers like
\begin{displaymath}
  {\cal Q}( {\cal Q}( \mbox{googol, googol} ) ,
            {\cal Q}( \mbox{googol, googol} ))
  \,\,\, ?
\end{displaymath}
Despite how much it breaks one's brain to think about it, the
number described by the above expression
{\em is still finite.}

\citet{Rayo19} discussed other possible generalizations (i.e.,
replacing ``first-order set theory'' by hypothetical higher-order
theories) and concluded by musing that
``Our quest to find larger and larger numbers has now morphed into
a quest to find more and more powerful languages!''

Another lesson that can be learned at this point is that, sometimes,
one needs to break out of a pattern rather than
just continuing it in the usual way.
In other words, we are about to see many examples of
``thinking outside the box,'' ``cutting the Gordian knot,'' or
even ``defeating the Kobayashi Maru.''

\section{Semi-Infinity: A Brief Rest Stop}
\label{sec:powertower}

Before we move on to the ``real deal,'' there is another useful
concept to examine that is both somewhat finite and somewhat infinite.
Consider the {\em infinite power tower} of exponentiation, which
can be expressed as a function $y$ of a positive real number $x$:
\begin{equation}   
  y(x) \,\, = \,\,
  x^{x^{x^{x^{\udots}}}} \,\, = \,\,
  x \, \uparrow\uparrow \, \infty \,\,\, ,
  \label{eq:tower1}
\end{equation} 
where this is probably the only time that we'll seriously use the
symbol $\infty$.
Here, we will consider an infinite number of exponentiations,
but we will assume only finite values for $x$.
Much more information about this function has been
given by, e.g., \citet{Kn81}, \citet{Mo19}, and \citet{Gd20}.
There are some values of $x$ for which $y(x)$ converges to a
constant and finite value, some values of $x$ for which $y(x)$
diverges to infinity, and some values of $x$ for which $y(x)$
keeps oscillating forever between two different solutions.

Another way to think about this is as an algorithm that requires
an infinite number of iterations:
\begin{equation}
  \left\{ \,\, \begin{array}{rcl}
    y_1 & \, = \, & x \\
    y_2 & \, = \, & x^x \\
    y_3 & \, = \, & x^{x^x}
  \end{array} \right\} 
  \,\,\,\,\,\,
  \mbox{and, in general,} 
  \,\,\,\,\,\,
    y_{n+1} \, = \, x^{y_n} \,\,\, ,
  \,\,\,\,\,\,
  \mbox{so that} 
  \,\,\,\,\,\,
  y(x) \, = \, \lim_{n \rightarrow \infty} \, y_n  \,\,\, .
\end{equation}
However, if we are examining a case in which $y(x)$ converges to a
finite value, then it means we ultimately reach a {\em fixed point}
in the iteration cycle: i.e., the ratio $y_{n+1}/y_n$ eventually
approaches 1.
Thus, there exists a self-consistent solution to the equation
\begin{equation}
  y \,\, = \,\, x^y \,\,\, .
  \label{eq:yxy_fixed}
\end{equation}
Notice that the $y$ on the right-hand side of Equation~(\ref{eq:yxy_fixed})
is already equal to an infinite chain of $x^{x^{\udots}}$, so adding one
more factor of ``$x$~to~the...'' onto the stack essentially changes nothing.
Thus, the whole thing remains equal to $y$ itself.
We will see many more examples of fixed-point relations, so it is
important to dig a little deeper into how they work.

\begin{figure}[!t]
\epsscale{1.05}
\plotone{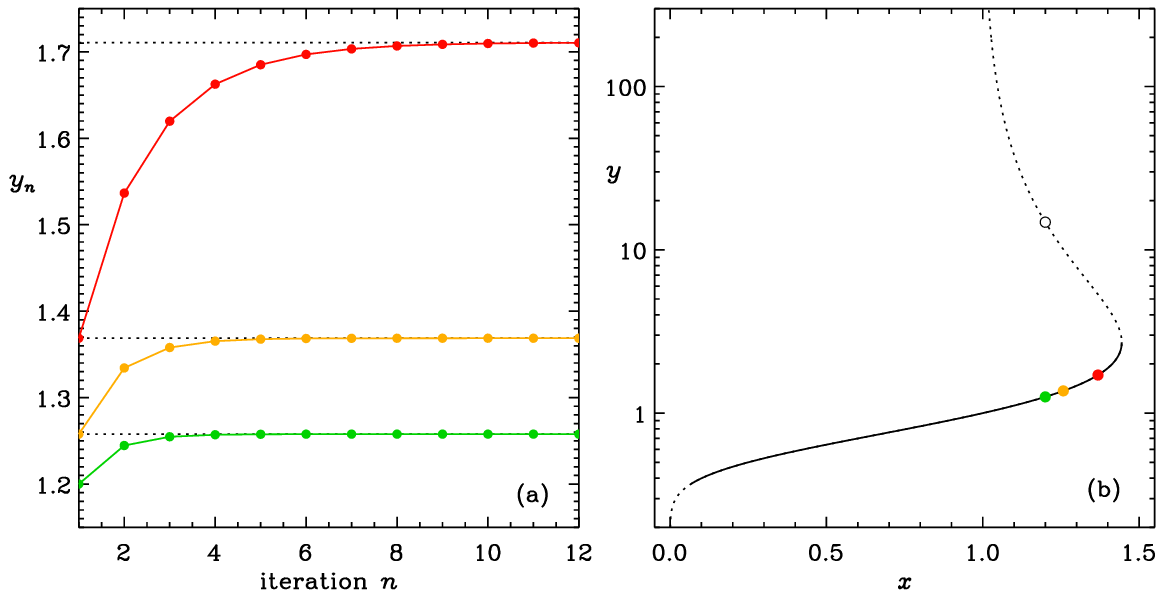}
\caption{(a) Example iterations $y_n$ from three different initial
values of $x$.
(b) Multi-valued fixed-point solution $y(x)$, with the stable solutions
shown as solid black line, unstable solutions as dotted lines, and
the three examples from panel~(a) shown with colored points.
\label{fig:tower}}
\end{figure}

We can start by noting that the range of $x$ values for which $y(x)$
converges to a finite value is given by
\begin{displaymath}
  e^{-e} \,\, < \,\, x \,\, < \,\, e^{1/e}
  \,\,\, , \,\,\,\,\,\,\,\,\,
  \mbox{i.e.,}
  \,\,\,\,\,\,\,\,\,
  0.065988 \,\, < \,\, x \,\, < \,\, 1.44467 \,\,\, .
\end{displaymath}
Thus if we choose a value within this range, we can see what happens when
we apply---and then repeat---the power-tower operation:
\begin{equation}
  \begin{array}{rcl}
    \mbox{For $x \, = \, 1.2$ ,}
    & y \, = \, x^{x^{\udots}} & \, = \, 1.257734... \\
    \mbox{For $y \, = \, 1.257734...$ ,}
    & z \, = \, y^{y^{\udots}} & \, = \, 1.368696... \\
    \mbox{For $z \, = \, 1.368696...$ ,}
    & t \, = \, z^{z^{\udots}} & \, = \, 1.710757...
  \end{array}
\end{equation}
Figure~\ref{fig:tower}(a) shows how these three examples converge
to their respective fixed-point solutions after a handful of iterations.

It may seem strange that the three example numbers ($y$, $z$, and $t$)
are not equal to one another.
After all, if $z$ is an infinite tower of $y$'s, and
$y$ is an infinite tower of $x$'s, shouldn't $z$ also be given by
the same infinite tower of $x$'s, as well?
The answer to that is no, due to the fact that exponentiation is
not associative.
In other words,
\begin{displaymath}
  y \,\,\, \mbox{is} \,\,\,\,
    x \expon
  ( x \expon
  ( x \expon
  ( x \expon
  ( x \expon
  ( x \expon
  ( x \expon
  ( x \expon
  ( x \expon \cdots ) ) ) ) ) ) ) )
  \,\,\,\,\,\,\,\,\,
  \,\,\,\,\,\,\,\,\,
\end{displaymath}
but
\begin{displaymath}
  z \,\, \mbox{is} \,\,
        ( x \expon ( x \expon ( x \expon \cdots ) ) ) \expon
  \Big\{ ( x \expon ( x \expon ( x \expon \cdots ) ) ) \expon
  \Big[ ( x \expon ( x \expon ( x \expon \cdots ) ) ) \expon \cdots
  \Big] \, \Big\}
\end{displaymath}
so $z \neq y$.

There is one additional thing to note about our first example of $x = 1.2$.
There is another value of $y$ that satisfies the fixed-point
condition of Equation~(\ref{eq:yxy_fixed}).
Try plugging in both $x=1.2$ and $y = 14.767453$ for yourself,
and you can show that $y= x^y$.
One can find this value by taking the fixed-point relation and
inverting it to solve for $x = y^{1/y}$.
Figure~\ref{fig:tower}(b) shows the result, and this ``hidden''
solution is illustrated with a white circle symbol.
Note that this new value of $y$ is most definitely {\em not} a solution
to the original power-tower condition of Equation~(\ref{eq:tower1}).
\citet{Mo19} explores many more interesting implications of these
multiple solutions, including their relationships to stability
and instability in dynamical systems.

We went through these examples to demonstrate two general facts about
fixed-point relations that will be useful to recall later:
\begin{itemize}
\item
If we have a fixed-point condition, we can generate
multiple examples of solutions by taking the output of a previous
solution as the input to the next.
\item
Sometimes a fixed-point condition ``takes on a life of its own'' and
exhibits solutions that never would have been apparent from the
original infinitely-iterated condition that gave rise to it.
\end{itemize}

\section{Countable Infinities}
\label{sec:countable}

Now we need to back up a bit and think about numbers in slightly different
ways.
It's helpful to think about numbers---and groups of numbers---as
{\em sets.}
We start with the set of all ``counting numbers,'' which we can define
as zero and all of the positive integers.
These will also be called the ``natural numbers,'' but that is
imprecise because some mathematicians exclude zero from that definition.
The set containing all of the natural numbers is often called $\mathbb N$,
and \citet{VN23} devised a particularly elegant way of representing
each member of $\mathbb N$ as sets unto themselves:
\begin{equation}
  \begin{array}{rll}
    0 \, = \, & \{ \, \} & \, = \,\, \emptyset \\
    1 \, = \, & \{ 0 \} & \, = \, \{ \emptyset \} \\
    2 \, = \, & \{ 0,1 \} & \, = \, \{ \emptyset , \{ \emptyset \} \} \\
    3 \, = \, & \{ 0,1,2 \} & \, = \, \{ \emptyset , \{ \emptyset \} ,
    \{ \emptyset , \{ \emptyset \} \} \} \\
    4 \, = \, & \{ 0,1,2,3 \} & \, = \, \{ \emptyset , \{ \emptyset \} ,
    \{ \emptyset , \{ \emptyset \} \} , 
    \{ \emptyset , \{ \emptyset \} , \{ \emptyset , \{ \emptyset \} \}\}\}
  \end{array}
  \label{eq:vonNeumann}
\end{equation}
and so on, where $\emptyset$ is the traditional symbol for the empty set.
The set for natural number $n$ contains $n$ members, and these members
correspond to the natural numbers from 0 to $n-1$.
Interestingly, this whole edifice is built upon
rearrangements of rearrangements of the empty set!
This kind of number theory has spurred physicists to think about our
actual universe as being constructed from nothing but ``bits'' of
pure information \citep[see, e.g.,][]{Wh90,Wi99}.

\subsection{Taking the Limit}
\label{sec:countable:lim1}

Let us ask the question:
What happens when we just start counting, and never stop?
Essentially, we start with zero, and then we repeatedly implement
operations of the command ``add one.''
We repeat, and repeat, and repeat, and
{\em then} we slyly say that we are skipping to the end
(i.e., taking the limit):
\begin{equation}
  0 \, , \, 1 \, , \, 2 \, , \, 3 \, , \,
  4 \, , \, 5 \, , \, 6 \, , \, \ldots \, , \, \omega \,\,\, .
  \label{eq:natural}
\end{equation}
\citet{C83} chose the symbol $\omega$ for the ``end'' of this sequence,
i.e., the number that follows all of the other finite numbers.
It is also called the first infinite {\em ordinal} number,
where we expand on the familiar idea of ordinal numbers as
those that tell us where numbers fall in an ordered sequence
(i.e., first, second, third, and so on).

Using the language of set theory, a compact way of defining this
first infinite ordinal is
\begin{equation}
  \omega \,\, = \,\, \bigcup_{n \in {\mathbb N}} \, n \,\,\,\, ,
\end{equation}
where the symbol $\in$ means ``is a member of.''
In other words, $\omega$ is the union (i.e., the full collection
of all members) of the sets corresponding to each of the finite
ordinals $n$ that make up the set of {\em all} natural numbers $\mathbb N$.

\begin{figure}[!t]
\epsscale{1.11}
\plotone{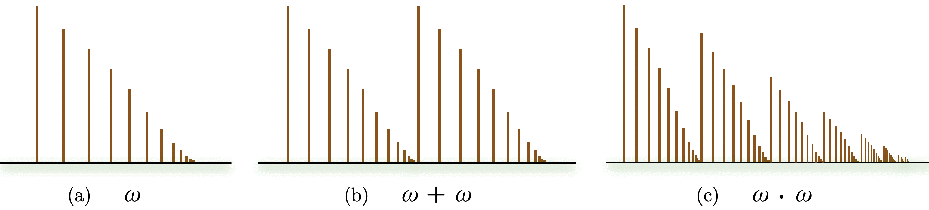}
\caption{Three infinite ordinals illustrated as telephone poles
receding toward the horizon.
\label{fig:poles}}
\end{figure}

Does $\omega$ really exist?
It does not have a known numerical value---which is a feature that
one might have assumed to be a key prerequisite for ``being a number.''
In a way, thinking about $\omega$ must involve a cognitive process
like chunking \citep{Gb01} or some degree of holistic or
Gestalt thinking \citep{Sm88}.
The human mind naturally assembles together related items into groups
or categories.
When we want to understand something, we often benefit from perceiving it
as a unitary whole, in addition to understanding that it consists of
smaller parts.
Thus, people have developed analogies for visualizing an infinite
number of steps in ways that do not require us to think about the
truly infinite sequence of numbers that defines it.
For example,
\begin{itemize}
\item
\citet{Ru82} discussed the analogy of ``speed-ups,'' in which we are
taking steps that get progressively smaller.
If the $n^{\rm th}$ step is of length $1/2^n$, then the total
distance traveled is finite, since
\begin{equation}
  1 \, + \, \frac{1}{2} \, + \,
  \frac{1}{4} \, + \,
  \frac{1}{8} \, + \,
  \frac{1}{16} \, + \, 
  \frac{1}{32} \, + \, \cdots \,\,\, = \,\,\, 2  \,\,\, .
  \label{eq:speedup}
\end{equation}
\item
Alternately, \citet{CG96} and \citet{Sh14} visualized these steps as
if one were standing on a long, straight road that has a row of
telephone poles alongside it (see Figure~\ref{fig:poles}(a)).
As one looks off into the distance, the poles recede toward the horizon,
and the entire collection can be perceived essentially as a finite thing.
\item
\citet{Th54} coined the phrase ``supertasks'' for infinitely iterated
steps that may or may not be possible to compress into a finite amount
of time or space.
Philosophers of science and mathematics have associated supertasks
with Zeno's paradox and other thought experiments that appear to
allow for mutually contradictory conclusions \citep[see also][]{Cl12}.

\end{itemize}

There is another basic property, which all natural numbers have,
but $\omega$ doesn't have: $\omega$ has no immediate predecessor.
All finite natural numbers can be called {\em successor ordinals}
because each number $n$ is the immediate successor of some other
known number $n-1$.
If we try to identify some natural number, which is both less
than $\omega$ and is its immediate predecessor, we fail.
For whatever natural number $n$ we identify, it has an immediate
successor $n+1$ which is not equal to $\omega$.
Thus, $\omega$ must be an example of a completely new type of ordinal
that is not a successor ordinal.
We now introduce a new category called {\em limit ordinals}
to describe $\omega$.
All ordinal numbers---finite or infinite---are either successor ordinals
or limit ordinals.
(Strictly speaking, zero is in a class all by itself.)

Where do we go from here?  Our impulse is always to keep counting.
One answer to the question ``what's bigger than infinity?''\  is
obviously ``infinity plus one!''
Despite that being a facetious answer, we can consider
continuing the original series of natural numbers with
\begin{equation}
  0 \, , \, 1 \, , \, 2 \, , \, 3 \, , \,
  4 \, , \, 5 \, , \, 6 \, , \, \ldots \, , \, \omega \, , \,
  \omega+1 \, , \,
  \omega+2 \, , \,
  \omega+3 \, , \,
  \omega+4 \, , \,
  \omega+5 \, , \, \ldots
  \label{eq:omplusn}
\end{equation}
Doing something like this may seem far more esoteric than merely
postulating the existence of $\omega$.
However, this continuation into additional ``transfinite'' numbers
does make sense, and we will discuss how to visualize and think
about them.

One thing we can note is that $\omega+1$ (and indeed any
$\omega + n$, for finite $n$) may be infinite, but it {\em does} have
an immediate predecessor.
Thus, now we have an example of an infinite successor ordinal.
Not all infinite numbers are limit ordinals.

At this point, many books divert into something a bit more
esoteric, which I don't think helps us that much.
This is the point that ``adding'' transfinite numbers is no longer
commutative,
\begin{equation}
  \mbox{i.e.,} \,\,\,\,\,\,\,\,\, 1 + \omega \, \neq \, \omega + 1
  \,\,\, .
\end{equation}
Yes, in one sense, it's true.  $1 + \omega$ starts with an extra finite
step, but ends in the same limiting place as just plain $\omega$.
On the other hand, $\omega + 1$ goes to that limit, and then takes
one {\em more} step.
These are ordinal numbers, and for them, ordering matters.
However, I'm going to now claim that, in another sense, this distinction
is a bit useless, and it can be shown that one can construct a
one-to-one correspondence between the members of the set
with $\omega$ members, the set with $1 + \omega$ members, and the
set with $\omega + 1$ members.
In that sense, these three sets are all of the same ``size.''

An analogy that is often employed to help us understand concepts
like $\omega + 1$ is {\em Hilbert's Hotel,}
a fanciful building with an infinite number of rooms, each numbered
according to the natural numbers: 0, 1, 2, 3, and so on
\citep[see, e.g.,][]{ES13,Kr14}.
Suppose, one day, the entire hotel is filled with guests (i.e., a
number of guests equal to $\omega$).
We can make a list of all guests and match them up to room numbers,
such that each guest has a unique room number, and no room is empty.
However, what happens when one more guest arrives?
Is there room for the newcomer?

\begin{figure}[!t]
\epsscale{1.01}
\plotone{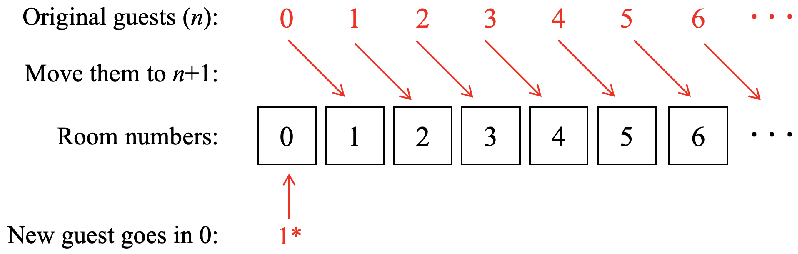}
\caption{Adding one guest to Hilbert's Hotel.
\label{fig:hilbert1}}
\end{figure}

Of course!
Figure~\ref{fig:hilbert1} shows how the hotel manager moves the guest in
Room~0 into Room~1, while simultaneously
moving the guest in Room~1 into Room~2,
the guest in Room~2 into Room~3, and so on.
Then, the new guest can have Room~0.
This means that if we made a new list of all of the guests, including
the new one, we are {\em still} able to match them up to room numbers that
correspond to the set of all natural numbers.
This is what we meant above by constructing a one-to-one correspondence
between $\omega$ and $\omega + 1$.
To paraphrase \citet{St16}, $\omega + 1$ is not actually larger
than $\omega$; it just is the number that comes next, in the ordered
list of ordinal numbers, immediately {\em after} $\omega$.

What, then, about the difference between $\omega + 1$ and $1 + \omega$?
\citet{Ch17} made the insightful analogy that Figure~\ref{fig:hilbert1}
can be interpreted in two different ways.
First, as we discussed above, is the case of starting with $\omega$ guests,
then encountering one more that wants to check in.
This requires all earlier $\omega$ guests to move.
Second, however, consider the situation where we started with that one
guest in room 0, and {\em then} $\omega$ guests arrived slightly later.
We could immediately place them into room numbers $n+1$ (for each
new guest numbered $n$), with nobody needing to move.
We obtained an outcome of the ``same size'' (i.e., a filled hotel)
each time, but the ordering of events was different.

\subsection{Taking More Limits}
\label{sec:countable:lim2}

At this point, I suggest going back and taking another look at
Equation~(\ref{eq:natural}), and then at Equation~(\ref{eq:omplusn}).
Just like we took the limit of Equation~(\ref{eq:natural}) to
get to $\omega$, we can take the limit of Equation~(\ref{eq:omplusn})
to obtain the quantity
\begin{displaymath}
  \omega \, + \, \omega \,\,\, .
\end{displaymath}
This is now a second example of a limit ordinal, and it is illustrated
using Conway \& Guy's telephone poles in Figure~\ref{fig:poles}(b).
The compact way to express this number in set-theory notation is
\begin{equation}
  \omega + \omega \,\, = \,\, \bigcup_{n \in {\mathbb N}} \,
  ( \omega + n ) \,\,\,\, .
\end{equation}

\begin{figure}[!t]
\epsscale{1.01}
\plotone{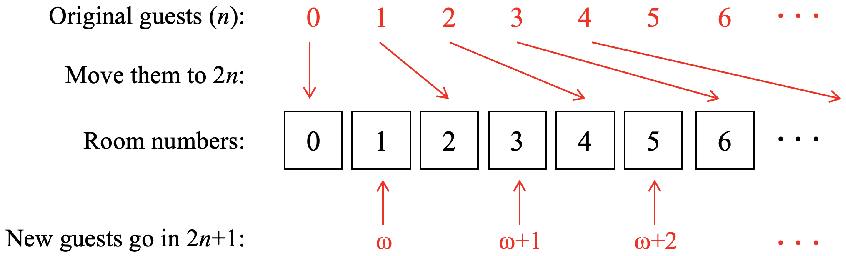}
\caption{Adding $\omega$ guests to Hilbert's Hotel.
\label{fig:hilbert2}}
\end{figure}

What do we mean by this?
Consider Hilbert's Hotel again.
As before, we start with it being full (with $\omega$ guests), but
now we want to allow {\em another} $\omega$ guests to stay there, too.
As shown in Figure~\ref{fig:hilbert2},
we can accomplish this by moving all existing guests (in rooms $n$)
to even-numbered rooms $2n$.
This frees up all of the odd-numbered rooms ($2n+1$) to be filled in
by the new guests.

Hold on for just a second.
This means that the set of all even or odd numbers,
despite being ``half'' of the set of all natural numbers, can be put
in a one-to-one correspondence with the full set of natural numbers.
This bit of weirdness was first pointed out by the philosopher
John Duns Scotus in the 1300s, so we may need to stew on it for a bit,
as well.
For now, let us be happy that we found a way to better visualize
$\omega + \omega$, and to show that it has a one-to-one correspondence
with $\omega$ itself.

We should note that mathematicians use a range of terms to describe
any two sets that have this state of one-to-one correspondence:
\begin{itemize}
\item
the mapping of elements is an exact {\em bijection}
(i.e., each element of the first set is paired with exactly one
element of the second set,
and each element of the second set is paired with exactly one
element of the first set),
which means there are no unpaired elements,
\item
they are equinumerous (they have essentially the same number of members),
\item
they are equipollent (they have the same ``power''),
\item
they have the same {\em cardinality.}
\end{itemize}
We will treat these terms as synonyms.

This same property applies to the set of all integers
(i.e., both negative and positive whole numbers).
Let us consider guests arriving to Hilbert's Hotel that are numbered
by integers $n$.
It's true that we cannot assign rooms to guests in a strictly increasing
order, as before, since the ``smallest'' integer corresponds to
$n \rightarrow -\infty$.
However, there is a scheme that can fill up the rooms in an orderly way.
We can start by putting guest $n=0$ into Room~0.
Then, guests with $n < 0$ can be given Rooms~$(-2n)$, and guests with
$n > 0$ can be given Rooms~$(2n-1)$:
\begin{displaymath}
 \begin{array}{rrrrrrrrrrr}
  \mbox{Rooms:}  \,\,\, &
  \, 0 &  1 &  2 &  3 &  4 &  5 &  6 &  7 &  8 & \,\, \cdots \\
  \mbox{Guests:} \,\,\, &
  \, 0 & +1 & -1 & +2 & -2 & +3 & -3 & +4 & -4 & \,\, \cdots
 \end{array}
\end{displaymath}
Each unique integer corresponds to a given unique natural number,
and each unique natural number corresponds to a given unique integer.
Thus, if we're ever asked ``how many integers are there?''\  we can
answer confidently: $\omega + \omega$.

After seeing many examples like this, we must realize that there
is something counterintuitive about infinite sets that isn't true
for finite sets:
an infinite set can be put into one-to-one correspondence with a
subset of itself!
In fact, it was only after mathematicians formally ``allowed'' this
to occur (using a formal axiom, which we will discuss more later)
that the theory of infinite numbers could truly proceed with a
firm logical foundation.

But we can go farther.
Note that our earlier limit of $\omega + \omega$ can be thought of as
$\omega \cdot 2$.
(Some books tidy up the notation and call it $2\omega$.
Although we're not being too pedantic about the non-commutativity
of infinite addition, I still want to stick with the most correct
version, which is $\omega \cdot 2$.)
Then, we keep adding by ones, as before, and we get to
\begin{equation}
  \omega \cdot 2 \,\, , \,\,
  \omega \cdot 2 + 1 \,\, , \,\,
  \omega \cdot 2 + 2 \,\, , \,\,
  \omega \cdot 2 + 3 \,\, , \,\,
  \ldots \,\, , \,\,
  \omega \cdot 2 + \omega  \,\,\, = \,\,\, \omega \cdot 3  \,\,\, .
\end{equation}
If we repeat this process again and again, we can get to
$\omega \cdot 4$, then $\omega \cdot 5$, and eventually take
an ``outer'' limit (sort of an umbrella ``limit of limits'')
to get $\omega \cdot \omega$.
We call this $\omega^2$, and we illustrate it in Figure~\ref{fig:poles}(c).
This is an infinite ordinal that corresponds to $\omega$ groups with
each group containing $\omega$ members.

We can show that this quantity applies to the positive rational
numbers (i.e., fractions $p/q$ with both $p$ and $q$ being
positive whole numbers).
The fact that a rational number is represented by a {\em pair}
of natural numbers seems to map quite nicely onto the concept of
$\omega^2$, which is what one gets after taking all permutations
of $\omega$ with $\omega$.
However, the rational numbers are still ``countable,'' and they have
the same cardinality as $\omega$ itself.
To show this, we need to show that {\em all} rationals can be enumerated
in an ordered list.

\begin{figure}[!t]
\epsscale{0.58}
\plotone{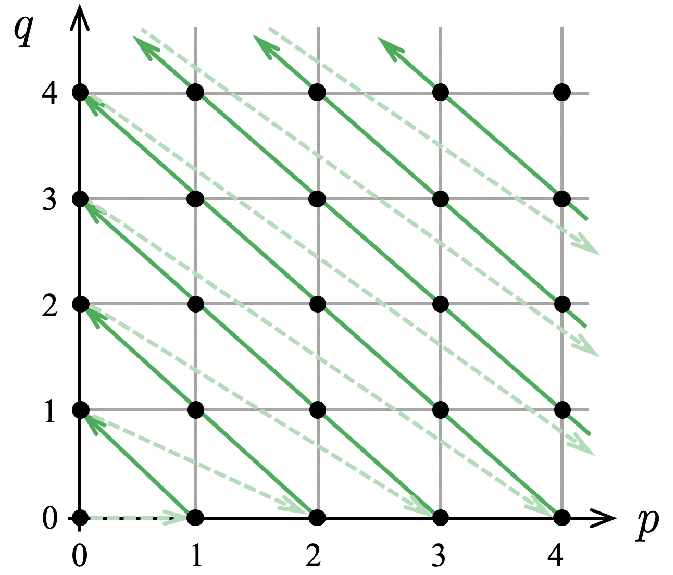}
\caption{One way to line up the rational numbers in a countable order.
\label{fig:rational}}
\end{figure}

The most common way to do this is to graph all possible numerators $p$
versus denominators $q$, as is shown in Figure~\ref{fig:rational}.
Then, we can start at the origin and trace diagonal lines through
successive groups with identical sums $p+q$.
As these lines pass through each unique pair ($p,q$), we assign new
natural numbers to them, in order.
Many books take pains to skip over fractions that are equivalent to
ones encountered earlier (i.e., one need not count 2/4 or 3/6 after
first encountering 1/2).
They also avoid rationals of the form $p/0$, $0/q$, and $0/0$.
However, I assert that it is fine to include all of these, since
there is still plenty of room for all of the unique, well-behaved
rational fractions.
The assignment of guests to rooms in Hilbert's Hotel is then:
\begin{displaymath}
 \begin{array}{rcccccccccccccccc}
  \mbox{Rooms:}  \,\,\, &
  \, 0 &  1 &  2 &  3 &  4 &  5 &  6 &  7 &  8 &
     9 & 10 & 11 & 12 & 13 & 14 &
     \,\, \cdots \\
  \mbox{Guests:} \,\,\, &
  \, 0/0 & 1/0 & 0/1 & 2/0 & 1/1 & 0/2 & 3/0 & 2/1 & 1/2 &
     0/3 & 4/0 & 3/1 & 2/2 & 1/3 & 0/4 &
     \,\, \cdots
 \end{array}
\end{displaymath}

Again, if we don't mind multiple-counting the duplicates and including
rationals with either $p=0$ or $q=0$, we can use
Cantor's pairing function
\begin{equation}
  n(p,q) \,\, = \,\, \frac{1}{2}(p+q)(p+q+1) + q  \,\,\,\, ,
\end{equation}
which generates the full set of natural numbers, in order, from
the traced sequence of points shown in Figure~\ref{fig:rational}.
The inverse function, which gives $p$ and $q$ for any $n$, is
\begin{equation}
  w \, = \, \left\lfloor \frac{\sqrt{8n+1}-1}{2} \right\rfloor
  \,\,\, , \,\,\,\,\,\,\,
  q \, = \, n - \frac{w^2 + w}{2}
  \,\,\, , \,\,\,\,\,\,\,
  p \, = \, w - n + \frac{w^2 + w}{2}
\end{equation}
where $\lfloor x \rfloor$ is the ``floor function'' of $x$,
or the greatest integer less than or equal to $x$.
If you have never seen these formulas before, I encourage you to
program them into a computer and
verify for yourself that they do what I claim they do.

Anyway, this one-to-one correspondence may seem more unusual than
the earlier ones, since there is an
``infinity's worth'' of rational fractions in between every pair of
integers on the number line.
In fact, even if you choose two rational fractions that are
very close to one another, but not equal, you can {\em still} squeeze
in an infinity's worth of other rational fractions between them!

Before we move on to even larger ordinals, we can note one more
interesting thing about $\omega^2$.
Because this is the same quantity as $\omega \cdot \omega$,
we can go back to Equation~(\ref{eq:multdef}), the definition of
multiplication, and write $\omega^2$ in another way:
\begin{equation}
  \omega^2 \,\,\, = \,\,\,
  \underbrace{\omega \, + \, \omega \, + \, \omega \, + \,
  \cdots}_{\omega \,\, {\rm times}}  \,\,\, .
\end{equation}
Because this is now an infinite chain of additions,
we can say that $\omega^2$ is the first ordinal $a$ for which
$\omega + a = a$.
In other words, putting another $\omega +$ in front
doesn't really change anything.
This is another example of a {\em fixed point,} a value
that doesn't change when subjected to a given transformation.

We can keep going.
We managed to construct $\omega^2$, so that means we can go on to
construct $\omega^3$, then $\omega^4$, and so on to another limit.
If we go back to Equation~(\ref{eq:expdef}), the definition of
exponentiation as repeated multiplication, it allows us to define
$\omega^{\omega}$ as the first ordinal $a$ for which
$\omega \cdot a = a$.
Yet another new flavor of fixed point!

The ordinal $\omega^{\omega}$ has a few interesting analogies.
\citet{Dv00} discussed the ``space'' of all possible
{\em polynomials} with natural-number coefficients and exponents.
There are $\omega$ possible zero-order polynomials of this kind
(i.e., constants),
and there are $\omega^2$ possible first-order polynomials (i.e., $ax+b$,
with the two constants $a$ and $b$ each varying over the full range of
natural numbers).
Then, there are $\omega^3$ possible second-order polynomials
($ax^2 + bx + c$), $\omega^4$ possible third-order polynomials,
and so on.
It's possible to say that the total number of polynomials of this kind is
\begin{equation}
  \omega \, + \, \omega^2 \, + \, \omega^3 \, + \, \omega^4
  \, + \, \cdots
  \,\,\, = \,\,\, \omega^{\omega}  \,\,\, ,
\end{equation}
but we haven't really dealt with the kind of weird
{\em ordinal addition} on display here.
I think it's possible to understand it, more or less, just as a limit.
If we continue the sequence described above, we see that the ``largest''
polynomial would be an $\omega^{\rm th}$ order polynomial.
When all of its coefficients are counted, there would be
$\omega^{\omega + 1}$ of these polynomials.
This number is equivalent to $\omega \cdot \omega^{\omega}$, and
this is kind of ``just another iteration'' of the fixed-point
transformation discussed above (i.e., $\omega \cdot a = a$).
Thus, we can think of $\omega \cdot \omega^{\omega}$ as being
essentially the same number as $\omega^{\omega}$.

Another interesting way to imagine $\omega^{\omega}$ is as an infinite
regress of countably infinite objects that encompass all prior ones.
Quoting \citet{Bz16},
\begin{quotation}
\noindent
First, imagine a book with $\omega$ pages. Then imagine an encyclopedia of
books like this, with $\omega$ volumes. Then imagine a bookcase containing
$\omega$ encyclopedias like this. Then imagine a room containing $\omega$
bookcases like this. Then imagine a floor
[of a]
library with $\omega$ rooms
like this. Then imagine a library with $\omega$ floors like this. Then
imagine a city with $\omega$ libraries like this.
And so on, {\em ad infinitum.}
\end{quotation}

We can now keep going and extend the sequence further.
What is the limit of the sequence
\begin{displaymath}
  \omega \,\,\, , \,\,\,
  \omega^{\omega} \,\,\, , \,\,\,
  \omega^{\omega^{\omega}} \,\,\, , \,\,\,
  \omega^{\omega^{\omega^{\omega}}} \,\,\, , \,\,\,
  \ldots \,\,\,\, ?
\end{displaymath}
Performing one more step of limit-taking, we can now think about
an infinite chain of $\omega$ to the power $\omega$ to the power
$\omega$ to the power... and so on.
This is a fixed-point defining an ordinal $a$ for which $\omega^a = a$.
We have seen this kind of power-tower before in
Section~\ref{sec:powertower}, but only for finite arguments.
Now we must confront its infinite cousin.
Cantor gave this ordinal the name $\varepsilon_0$, and we can define it as
\begin{equation}
  \varepsilon_0 \,\, = \,\,
  ^{\omega}\omega \,\, = \,\,
  \omega^{\omega^{\omega^{\omega^{\udots}}}}
\end{equation}
i.e., this is $\omega$ tetrated to the $\omega$, and it can be
expressed in set-theory notation as
\begin{equation}
  \varepsilon_0 \,\, = \,\,
  \bigcup_{n \in {\mathbb N}} \, ( \omega \uparrow \uparrow n ) \,\,\,\, .
\end{equation}

One can keep going with even higher analogues to the mathematical
operations, like
\begin{displaymath}
  ^{^{^{^{\ddots}\omega}\omega}\omega}\omega
\end{displaymath}
(i.e., $\omega$ pentated to the $\omega$), and so on,
but Cantor took a different route.
He noted that $\varepsilon_0$ was not the only solution to the
fixed-point equation $\omega^a = a$.
He showed that the ``next'' solution could be written as
\begin{equation}
  \varepsilon_1 \,\, = \,\,
    \varepsilon_{0}^{\varepsilon_{0}
  ^{\varepsilon_{0}^{\varepsilon_{0}^{\udots}}}}
\end{equation}
and, in some sense, $\varepsilon_1 > \varepsilon_0$.
A general series of $\varepsilon$ numbers can be written as
\begin{equation}
  \varepsilon_{m+1} \,\, = \,\,
  \bigcup_{n \in {\mathbb N}} \,
  ( \varepsilon_{m} \uparrow \uparrow n ) \,\,\, ,
  \,\,\,\,\,\,\,\,\,\,\,\,\,\,\,
  \mbox{for $m \geq 0$ .}
  \label{eq:epsilonm}
\end{equation}
Where does the infinite series of $\varepsilon$ numbers lead?
Well, one can imagine taking the limit to obtain $\varepsilon_{\omega}$,
and then climbing the ladder of infinite ordinals once again in the
subscript.
Eventually, we would be forced to contemplate
$\varepsilon_{\varepsilon_0}$.
Continuing this process would lead us to
\begin{equation}
  \varepsilon_{\varepsilon_{\varepsilon_0}}
  \,\,\, ,
  \,\,\,\,\,\,\,\,\,
  \mbox{then to}
  \,\,\,\,\,\,\,\,\,
  \varepsilon_{\varepsilon_{\varepsilon_{\varepsilon_0}}}
  \,\,\,\,\,\,\,\,\,
  \mbox{then, ultimately, to}
  \,\,\,\,\,\,\,\,\,
  \varepsilon_{\varepsilon_{\varepsilon_{\varepsilon_{\varepsilon_{\ddots}}}}}
\end{equation}
which some have named $\alpha$.
This new ordinal number $\alpha$
satisfies a new kind of fixed-point condition
($\varepsilon_{\alpha} = \alpha$) that doesn't involve iterated
repetition of basic arithmetic operations.
Instead, it depends
on the subtle meaning of the subscript in Equation~(\ref{eq:epsilonm})
and what it means to ``repeatedly iterate'' by applying it to itself.

One can stretch the brain to consider going even further, but it is
important to note that all of these quantities still have the
same cardinality as $\omega$, and they are still all {\em countable.}
All of these expansions---even up through the $\varepsilon$'s and the
$\alpha$'s---still involve finding new ways to squeeze ever more
guests into Hilbert's Hotel.
We have not yet figured out if there are any numbers truly bigger,
but we will soon see there exist numbers of guests so large that they
{\em cannot} fit into the hotel.
For now, though, we will use the symbol $\aleph_0$, often
pronounced ``aleph-null,'' to describe the specific level of
cardinality that we have been dealing with so far
\citep[see, e.g.,][]{Ft12}.

\section{The Real Continuum}
\label{sec:real}

We may have started out thinking that the rational numbers
``fill in the number line'' in an unfathomably dense way.
But, no, there are several strikes against that idea.
First, there's the notion that they are countable
(i.e., that they have a one-to-one correspondence with the natural
numbers).
This starts to reveal their limitations.
Then, of course, there are famous cases of {\em irrational} numbers
like $\pi$ and $\sqrt{2}$ for which it's been proven impossible to
write as rational fractions \citep[see, e.g.,][]{Fl06}.
Thus, it's the real numbers that truly fill in the number line,
and thus they're traditionally called {\em the continuum.}

Are the real numbers countable?
If we try to list them as ``merely'' their decimal expansions, then,
sure, they're countable.
In other words, one can come up with a way to arrange these decimal
expansions in an ordered list, and extend that list to ``infinity.''
But, in reality, {\em there are more real numbers than that,}
and thus the true set of {\em ALL} real numbers is uncountable.
We will offer up a proof of this provocative statement in a bit.

\begin{figure}[!t]
\epsscale{0.76}
\plotone{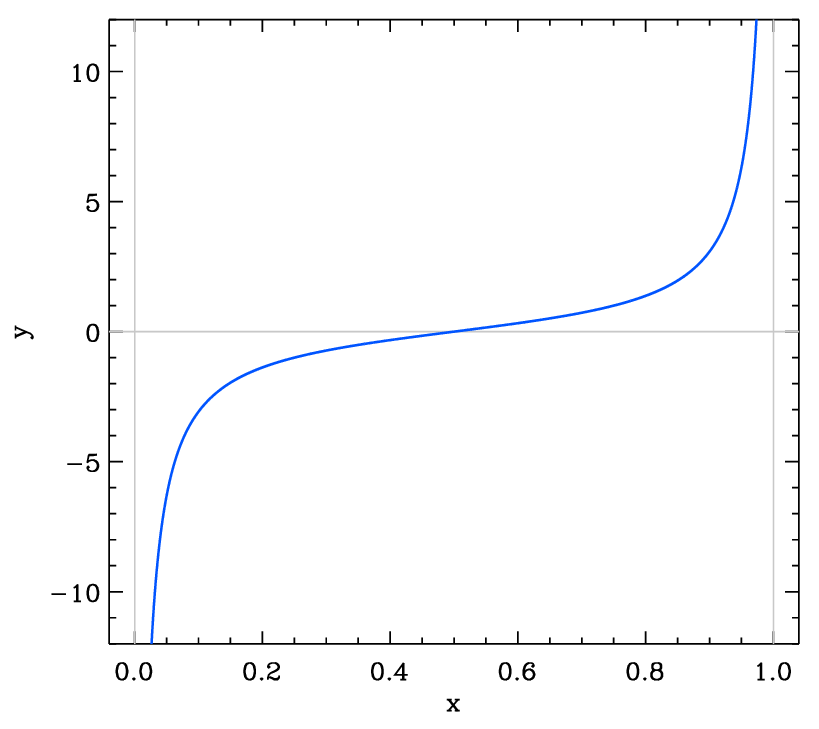}
\caption{A function that gives a one-to-one mapping of the reals
in the range $[0,1]$ to the full range of reals.
\label{fig:tanfunc}}
\end{figure}

First, though, we can make life easier for ourselves in one specific way.
We actually do not need to consider the {\em entire} set of real
numbers from $-\infty$ to $+\infty$.
We only need to deal with a finite ``line segment'' of the real
number line, which we choose as the space between 0 and 1.
Consider a real number $x$ that can take on any value between 0 and 1,
and compute the function
\begin{equation}
  y(x) \,\, = \,\, \tan \left[ \frac{\pi}{2} \Big( 2x - 1 \Big) \right]
  \,\,\, .
\end{equation}
Figure~\ref{fig:tanfunc} shows this function, and it is clear that
$y$ covers the full range of real numbers from $-\infty$ to $+\infty$.
Every real value of $x$ corresponds to a unique real value of $y$, and
every real value of $y$ corresponds to a unique real value of $x$.
In other words, this function gives an exact one-to-one mapping
between the reals in the range $[0,1]$ and the full set of reals.
Thus, everything that we prove for the interval $[0,1]$ will also
be true for the entire real line.\footnote{%
If this argument was not completely convincing, there are some
other good ones given by \citet{Ch17}.}

Our job now is to come up with some nice, ordered list that enumerates
all possible decimal expansions of the real numbers in this range.
Maybe we start with 0 (which is really 0.0000000...), and then
include the nine numbers describing tenths
($0.1, 0.2, 0.3, \ldots 0.9$), and then
the 99 numbers describing hundredths (0.01 to 0.99), then
the 999 numbers describing thousandths (0.001 to 0.999), and so on.
By gradually filling in all trailing zeros, it's essentially the same
as counting up the natural numbers like in Equation~(\ref{eq:natural}).

No matter what system we use to enumerate these decimals, we can
create a table that looks like:
\begin{displaymath}
  \begin{array}{cl}
    x_1 \,\, = \,\, & 0. d_{11} d_{12} d_{13} d_{14} d_{15} \,\, \ldots \\
    x_2 \,\, = \,\, & 0. d_{21} d_{22} d_{23} d_{24} d_{25} \,\, \ldots \\
    x_3 \,\, = \,\, & 0. d_{31} d_{32} d_{33} d_{34} d_{35} \,\, \ldots \\
    x_4 \,\, = \,\, & 0. d_{41} d_{42} d_{43} d_{44} d_{45} \,\, \ldots \\
    x_5 \,\, = \,\, & 0. d_{51} d_{52} d_{53} d_{54} d_{55} \,\, \ldots \\
    \vdots & \mbox{\hspace*{0.18in}}
             \vdots \mbox{\hspace*{0.19in}}
             \vdots \mbox{\hspace*{0.19in}}
             \vdots \mbox{\hspace*{0.19in}}
             \vdots \mbox{\hspace*{0.19in}}
             \vdots
  \end{array}
\end{displaymath}
and each digit $d_{ij}$ is known.
If this matrix extends forever in both directions (i.e., both going down
and going to the right) then surely it must contain every possible
real number, right?

Surprisingly, the answer is no.
In 1873, Georg Cantor showed that there are many possible real numbers
that are {\em not} contained in this enumerated list.
Here we will describe Cantor's {\em diagonal proof,} but this was only
one of several different proofs of this key idea that he developed.

Let us take the known digits $d_{ij}$ and define some new digits:
\begin{equation}
  D_k \,\, = \,\, \left\{
  \begin{array}{ll}
    d_{kk} + 1 \,\, , \,\, & \mbox{if $d_{kk} \neq 9$,} \\
    0 \,\, , \,\, & \mbox{if $d_{kk} = 9$ .} \\
  \end{array} \right.
\end{equation}
In other words, we take the digits that go down the diagonal in the
above table ($d_{11}$, $d_{22}$, $d_{33}$, $d_{44}$, and so on)
and create a new series of digits $D_k$ by adding 1 to each, and
cycling back to zero if the digit was 9.
Then, we create a new real number
\begin{displaymath}
  x_{\rm new} \,\, = \,\, 0.D_1 D_2 D_3 D_4 D_5 \,\, \ldots \,\,\,\, .
\end{displaymath}
Does this new real number occur in the table?
Well, we know that it differs from $x_1$ in the first digit after
the decimal, and that it differs from $x_2$ in the second digit after
the decimal, and that it differs from $x_3$ in the third digit after
the decimal, and so on.
The real number $x_{\rm new}$ cannot occur in our enumerated list!

Note that our methods of constructing both the ordered list
($x_1, x_2, x_3, \ldots$) and $x_{\rm new}$ were quite arbitrary.
There are literally an infinite number of other ways we could have
chosen to carry out this process, and $x_{\rm new}$ still would not
appear in the ordered list.

A key implication of this result is that there can never be a
countable enumeration of the real numbers that actually, exhaustively,
includes all of them.
In fact, it would seem that there are infinitely {\em more} of these
uncountable real numbers (all of them being irrational)
than there are countable ones.
If you were to randomly throw a dart at the real number line, the
probability that you would hit a rational number is so infinitesimally
small that it is indistinguishable from zero.
(See Appendix~\ref{sec:appA} for the outlines of a few proofs
that show this is the case.)

Using the symbol ${\mathfrak c}$ to refer to the cardinality (size)
of the real continuum, we have now shown that
\begin{equation}
  {\mathfrak c} \, > \, \aleph_0 \,\,\, .
\end{equation}
Amazingly, this demonstrates that there exist fundamentally different
``sizes of infinity.''

Cantor also showed that the multi-dimensional continuum
(i.e., the area bounded by any plane, the volume bounded by
any cube, and so on) has the same cardinality $\mathfrak c$
as the one-dimensional continuum of real numbers discussed above.
Trying to prove this was a long-time source of consternation to
Cantor, who said in a letter {\em ``je le vois, mais je ne le crois pas''}
(I see it, but I don't believe it).

However, we can start to understand this, for the
2D plane that covers the square region $[0,1]$ in $x$ and
$[0,1]$ in $y$, by finding a one-to-one correspondence between the
real numbers and the coordinates $(x,y)$.
Let us write out the decimal expansions of the coordinates,
\begin{displaymath}
  x \,\, = \,\, 0. a_1 a_2 a_3 a_4 a_5 a_6 \, \ldots
  \,\,\,\,\,\,\,\,\,\,\,\,
  \mbox{and}
  \,\,\,\,\,\,\,\,\,\,\,\,
  y \,\, = \,\, 0. b_1 b_2 b_3 b_4 b_5 b_6 \, \ldots
\end{displaymath}
and then interleave the digits, like when trying to fit $\omega$
additional guests into Hilbert's Hotel:
\begin{displaymath}
  r \,\, = \,\, 0. a_1 b_1 a_2 b_2 a_3 b_3 a_4 b_4 a_5 b_5 a_6 b_6
  \, \ldots
\end{displaymath}
Everything that you can do with $x$ and $y$---including feeding them
into Cantor's diagonal proof to generate uncountably
more numbers---can be done with $r$, as well.
Thus, every real number $r$ can be made to correspond to a unique
point in the 2D plane.
Another way to construct this kind of one-to-one mapping is to define
a fractal-like {\em space-filling curve} \citep{Sg94}
that, when iterated an infinite number of times, can eventually
wiggle around to reach every point in the 2D plane, while simultaneously
maintaining a known 1D path that is mappable to the real number line.

Lastly, we note that the set of all {\em complex numbers} is
essentially the same thing as the full set of points in a
two-dimensional plane.
Thus, it is also true that the complex numbers also have
cardinality $\mathfrak c$.

\section{Uncountable Infinities}
\label{sec:uncountable}

Are there other ways to construct a quantity with cardinality
greater than $\aleph_0$?
Trying to do this opens up the way to the ``higher'' infinities.

\subsection{Alephs and Beths}
\label{sec:uncountable:AB}

At this point, we are going to start talking more explicitly
about {\em cardinal numbers} such as $\aleph_0$ and $\mathfrak c$.
This is really the same concept as cardinality that we have been
discussing for a while, but we are now labeling each size-class with
a unique number.
Note that each cardinal number corresponds to a given lowest ordinal
of that size-class (i.e., $\aleph_0$ corresponds to $\omega$),
but there are multiple ordinal numbers that correspond to a given
cardinal (i.e., $\omega + 1$, $\omega \cdot \omega$, $\omega^{\omega}$,
and $\varepsilon_0$ all have cardinality $\aleph_0$, too).

Mathematicians define $\aleph_1$ as the ``next largest'' cardinal number
after $\aleph_0$, and $\aleph_2$ as the next largest after that, and so on.
At this point, you may be wondering how we can guarantee that successively
larger and larger cardinalities actually exist.
If so, stay tuned for just a few paragraphs.

Another way to conceptualize the $\aleph$ numbers is to define the
successor operation using superscript $+$.
For any finite natural number, $n^{+} = n+1$.
For an infinite cardinal $\kappa$, the next largest cardinal
is defined as $\kappa^{+}$, and this is {\em not} the same
thing as $\kappa + 1$.
Thus, the $\aleph$ notation gives us a way to define each successor as
\begin{equation}
  \aleph_{n+1} \, = \, \aleph_n^{+}  \,\,\, .
\end{equation}
Because $\aleph_0$ corresponds to a {\em countably} infinite
cardinality, we see that the succeeding $\aleph_n$ cardinals
(for $n \geq 1$) must correspond to sets that are
{\em uncountably} infinite.
Thus, $\aleph_1$ describes the smallest number for which it
is {\em impossible} to fit that number of guests into Hilbert's Hotel.

When comparing cardinalities, we see that
\begin{equation}
  \aleph_0 \, < \, \aleph_1 \, < \, \aleph_2 \, < \,
  \aleph_3 \, < \, \cdots
  \label{eq:alephchain}
\end{equation}
One way to think about any arbitrary $\aleph_n$ is that this is the
cardinal number for which there exist $n$ smaller $\aleph$-numbers.
Because there cannot be any gaps in the chain of $\aleph$-numbers,
\citet{Zu74} noted that {\em every} infinite cardinal
must be an $\aleph$-number.
Lastly, we can quote \citet{CG96} to describe the first few
aleph numbers as follows:
\begin{itemize}
\item
$\aleph_0$ is the number of finite ordinal numbers.
\item
$\aleph_1$ is the number of ordinal numbers that are either finite
or in the $\aleph_0$ class.
\item
$\aleph_2$ is the number of ordinal numbers that are either finite
or in the $\aleph_0$ or $\aleph_1$ classes.
\end{itemize}
and so on.

Just going by the above definitions, it is difficult to
understand how to visualize or construct anything that corresponds
to these higher aleph-numbers.
Also, as we hinted above, we still aren't quite sure if there's
an actual necessity for anything to truly {\em exist} at any given
higher-aleph cardinality.

\begin{figure}[!t]
\epsscale{0.98}
\plotone{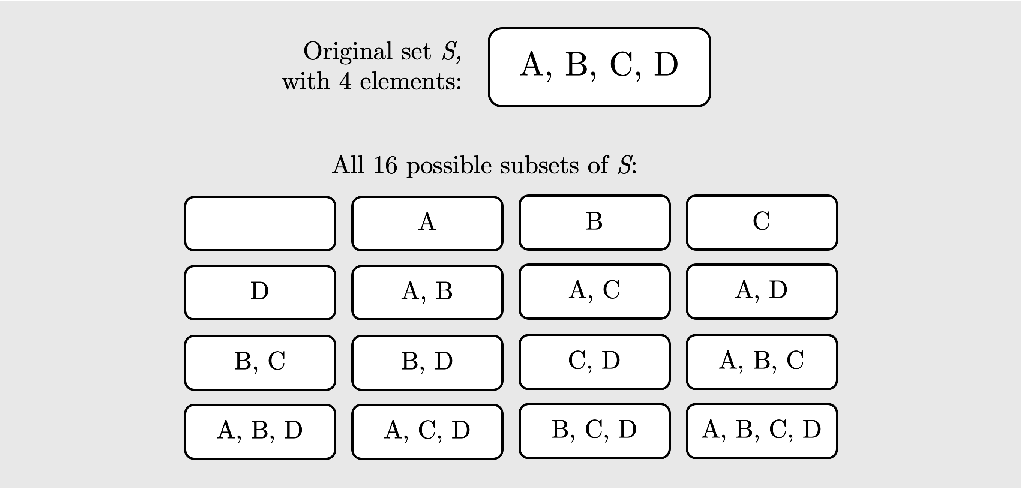}
\caption{Enumeration of all possible subsets of an example finite set.
\label{fig:powerset4}}
\end{figure}

A path forward---which guarantees that we can keep building new rungs
onto the top of this ladder as we climb it---is found by considering
the {\em Power-Set} of some set $S$.
This term is defined as the set of all possible subsets of $S$
(including the empty set and the ``full'' original set).
Figure~\ref{fig:powerset4} illustrates this with a finite set
that has 4 elements and $2^4 = 16$ possible subsets.
Note that membership matters, but not ordering.

Thus, as we move on to name sets by their size, we can state that
the cardinality of the power-set of $S$ is $2^S$, and it's
also true that $2^S > S$.
We will sometimes see the symbol ${\mathscr P}(S)$ used to express the
power-set of $S$.
The number $2^S$ occurs because there are $S$ independent Yes/No
decisions to make about whether
each member of the original set should be in a given subset.
Many authors, such as \citet{Ru82} and \citet{Ch17},
have illustrated this as a branching tree.

So, when going from finite sets to infinite sets, it's also true
that $2^{\aleph_0}$ (i.e., the power-set of $\aleph_0$) has a
greater cardinality than $\aleph_0$.
In other words, the power-set of the natural numbers isn't countable.
This assertion is known as {\em Cantor's theorem,} and we can prove it
with a diagonal argument very similar to the one discussed above.
We would like to enumerate all possible subsets of the set of all
natural numbers.
For each subset, we need to provide a unique set of Yes/No answers to
$\omega$ independent questions: is 0 in the subset?\  is 1 in
the subset?\  is 2 in the subset?\  and so on.
The table below shows one way of enumerating these subsets:
\begin{displaymath}
  \begin{array}{c|ccccccccc}
      & 0 & 1 & 2 & 3 & 4 & 5 & 6 & 7 & \cdots \\
    \hline
  S_0 & \mbox{N} & \mbox{N} & \mbox{N} & \mbox{N} &
        \mbox{N} & \mbox{N} & \mbox{N} & \mbox{N} & \cdots \\
  S_1 & \mbox{\bf Y} & \mbox{N} & \mbox{N} & \mbox{N} &
        \mbox{N} & \mbox{N} & \mbox{N} & \mbox{N} & \cdots \\
  S_2 & \mbox{N} & \mbox{\bf Y} & \mbox{N} & \mbox{N} &
        \mbox{N} & \mbox{N} & \mbox{N} & \mbox{N} & \cdots \\
  S_3 & \mbox{\bf Y} & \mbox{\bf Y} & \mbox{N} & \mbox{N} &
        \mbox{N} & \mbox{N} & \mbox{N} & \mbox{N} & \cdots \\

  S_4 & \mbox{N} & \mbox{N} & \mbox{\bf Y} & \mbox{N} &
        \mbox{N} & \mbox{N} & \mbox{N} & \mbox{N} & \cdots \\
  S_5 & \mbox{\bf Y} & \mbox{N} & \mbox{\bf Y} & \mbox{N} &
        \mbox{N} & \mbox{N} & \mbox{N} & \mbox{N} & \cdots \\
  S_6 & \mbox{N} & \mbox{\bf Y} & \mbox{\bf Y} & \mbox{N} &
        \mbox{N} & \mbox{N} & \mbox{N} & \mbox{N} & \cdots \\
  S_7 & \mbox{\bf Y} & \mbox{\bf Y} & \mbox{\bf Y} & \mbox{N} &
        \mbox{N} & \mbox{N} & \mbox{N} & \mbox{N} & \cdots \\
  \vdots & \vdots   & \vdots   & \vdots   & \vdots   &
        \vdots   & \vdots   & \vdots   & \vdots   & 
  \end{array}
\end{displaymath}
Then, like before, we can construct a new subset $S_{\rm new}$ out
of the diagonal members of the above matrix, after swapping Y for N
and N for Y.

A word of caution: In the particular example enumeration shown above,
I think $S_{\rm new}$ will consist of nothing but Y's.
I probably chose the most boring possible enumeration of Y and N,
and nearly all other books that discuss this proof give more
interesting enumerations.
However, for those, one usually has to take it on faith that there is
an actual {\em pattern} that eventually accounts for all countable
combinations.
For mine, at least, that pattern is evident.
(If that pattern isn't evident, I suggest looking up how to write
the natural numbers 0, 1, 2, 3, 4, $\ldots$ using base-two, or binary
digits instead of base-ten.)

Whatever $S_{\rm new}$ looks like, it cannot be a member of the
originally enumerated list.
Thus, the {\em complete} set of all subsets of the natural numbers
must be uncountable, and it must have a cardinality bigger
than $\aleph_0$.
We have now proven that
\begin{equation}
  2^{\aleph_0} \,\, > \,\, \aleph_0 \,\,\,\, .
\end{equation}

To help put these quantities in perspective,
some mathematicians define the ``beth'' numbers, such that
\begin{equation}
  \begin{array}{c}
    \beth_0 \, = \, \aleph_0 \\
    \beth_1 \, = \, 2^{\beth_0} \\
    \beth_{n+1} \, = \, 2^{\beth_n} \,\,\, , \,\,\,\,\,\,\,\,
    \mbox{and so on.}
  \end{array}
\end{equation}
In other words, each new $\beth$ number has the cardinality of the
power-set of the previous $\beth$ number.
Thus, when comparing cardinalities,
\begin{equation}
  \beth_0 \, < \, \beth_1 \, < \, \beth_2 \, < \,
  \beth_3 \, < \, \cdots
  \label{eq:bethchain}
\end{equation}
Unfortunately, though, it is not easy to figure out how
Equations~(\ref{eq:alephchain}) and (\ref{eq:bethchain})
relate to one another.
In other words, we still do not know how the $\aleph$ numbers
and the $\beth$ numbers are {\em interleaved} with one another.

\subsection{Back to the Continuum}
\label{sec:uncountable:C}

For now, it's helpful to go back to thinking about the
real continuum $\mathfrak c$.
Notice that the two diagonalization proofs given above told
us that ${\mathfrak c} > \aleph_0$ and that
$2^{\aleph_0} > \aleph_0$.
The details of those proofs were almost identical to one another.
In fact, they would have been {\em exactly} identical if we had
expressed the real numbers using binary digits, rather
than decimal digits!
Thus, because these two identical proofs essentially describe a
one-to-one correspondence, it is true that
\begin{equation}
  {\mathfrak c} \, = \, 2^{\aleph_0} \, = \, \beth_1 \,\,\, ,
  \label{eq:cbeth1}
\end{equation}
i.e., the cardinality of the real continuum is equal to the
cardinality of the power-set of the natural numbers.

Equation~(\ref{eq:cbeth1}) always kind of amazes me, but set theorists
seem to take it in stride.
As we will see below, there remain some deeply unsolvable mysteries
about the real continuum.
However, it is nice that we can hang our hats on
Equation~(\ref{eq:cbeth1}), which tethers the big, unfathomable
continuum to something specific and grounded about the
more-well-behaved natural numbers.

Georg Cantor strongly believed an additional thing to be true.
This thing is called {\em Cantor's continuum hypothesis} (CH):
\begin{equation}
  {\mathfrak c} \, = \, \aleph_1 \,\,\, .
\end{equation}
In other words, Cantor believed that there are no ``other''
cardinals squeezed in between $\aleph_0$ and ${\mathfrak c}$.
This is a tempting idea because there are deep conceptual
similarities between $\aleph_1$ and $2^{\aleph_0}$.
Because $\aleph_1$ is the cardinality (size) of the group of all
numbers that are either finite or in the $\aleph_0$ category,
it means that it's equivalent to the
{\em total number of countably infinite numbers.}
Similarly, because ${\mathfrak c} = 2^{\aleph_0}$ corresponds to the
power-set of $\aleph_0$, it seems to be equivalent to the
{\em total number of countable sets.}
Given the \citet{VN23} correspondence between numbers and sets
that we showed in Equation~(\ref{eq:vonNeumann}),
why would one ever doubt that ${\mathfrak c} = \aleph_1$?

In more recent years, there has arisen an even broader idea, called the
generalized continuum hypothesis (GCH), that proposes
\begin{equation}
  \aleph_{n+1} \, = \, 2^{\aleph_n} \, = \, \beth_{n+1}
  \,\,\,\,\,\,\,\,\,\,\,\,\,\,
  \mbox{for any finite $n$~.}
\end{equation}
In other words, GCH is the assertion that there are no
levels of cardinality intermediate between any $\aleph_n$ and
its corresponding power-set $2^{\aleph_n}$.
If GCH is true, then taking the immediate successor of a cardinal
is the {\em same thing} as taking the power-set of that cardinal.
This also makes sense because it is essentially what
Equation~(\ref{eq:vonNeumann}) does for all finite numbers.

Unfortunately, using the currently developed version of set theory
(Zermelo–Fraenkel theory, augmented with the Axiom of Choice, or ZFC),
there is no way to conclusively prove or disprove CH or GCH.
In other words, they're ``undecidable.''
Quoting \citet{Ru82},
\begin{quotation}
In 1940, Kurt G\"{o}del was able to show that CH is consistent with ZFC.
He showed that one can never prove ${\mathfrak c} \neq \aleph_1$ from
the axioms of ZFC. This does not mean that Cantor was {\em right}---it
only means that he was {\em not provably wrong} on the basis of ZFC.

In 1963, Paul Cohen proved that the negation of CH is consistent with ZFC.
He showed that one can never prove ${\mathfrak c} = \aleph_1$ from
the axioms of ZFC.
This does not mean that Cantor was {\em wrong}---it only means that we
{\em cannot prove that he was right,} using only the axioms of ZFC.

[...]$\,$
The situation is a little like asking what Scarlett O'Hara did after
the end of {\em Gone with the Wind...} one can consistently write a sequel
in which she gets back together with Rhett, and one can consistently write
a sequel in which she never sees Rhett again. But the book itself does not
give us enough information to draw either of these conclusions with
certainty. In the same sense, ZFC is not a complete enough description
of the universe of set theory to tell us what the power of the continuum is.
\end{quotation}
In the preface of the updated 2005 edition of \citet{Ru82}, there's
a discussion of recent work that may be starting to point
to the plausibility of ${\mathfrak c} = \aleph_2$
\citep[see also][]{For88,Bek91}.

Before we move on, we can say a few more words about ZFC.
It is much more powerful and internally consistent than the earlier
versions of set theory that developed at the end of the 19th century.
Like many other mathematical theories, it is built on a foundation
of {\em axioms.}
These are assertions, or postulates, that seem so self-evidently true
that they are just assumed at the outset and never actually proven.
ZFC uses several prosaic and sturdy axioms from finite set theory, and
it also includes the ``Axiom of Choice'' mentioned above---which took
a few decades for mathematicians to wrap their heads around---{\em{and}}
it adds one more:
\begin{quotation}
\noindent
{\em The Axiom of Infinity:} There exist sets that have a
one-to-one correspondence with subsets of itself.
\end{quotation}
(This is just one of many ways to state it.)
Without this axiom, one cannot rigorously work with infinite
concepts like $\omega$ and the $\aleph$ and $\beth$ numbers.

The process of tacking-on new axioms is something that we will
encounter a few more times before the end of our journey.
Strangely, this process is both constricting and freeing.
It is constricting because the more axioms one adopts, the more
rules one must obey, and this limits us to a smaller volume
of theoretical ``space'' in which to play.
However, it is freeing because these new axioms allow us to
confidently define and prove {\em new} things that would only be
uncertain and undecidable without them.

\section{Beyond the Continuum?}
\label{sec:beyondC}

Recall that ${\mathfrak c} = \beth_1$.
Thus, if the CH is true, then we may have $\beth_1 = \aleph_1$.
However, if the CH is not true, then the continuum must sit somewhere
``higher up'' in the chain of uncountable cardinals, such that
$\beth_1 > \aleph_1$.
This implies a more general inequality:
\begin{equation}
  \beth_n \,\, \geq \,\, \aleph_n  \,\,\,\,\,\,\,\,\,\,\,\,
  \mbox{for any $n$ .}
  \label{eq:bgeqa}
\end{equation}
Thus, if GCH is not true, then taking the power-set of a cardinal
(i.e., going up one step in the $\beth$ chain) represents a
fundamentally {\em bigger step} than just taking the immediate successor
of that cardinal (i.e., going up one step in the $\aleph$ chain).

Is there anything interesting to say about $\aleph_2$ or $\beth_2$?
Consider the power-set of the real numbers; i.e., the set of all
possible subsets of real numbers.  This has cardinality
\begin{equation}
  2^{\mathfrak c} \,\, = \,\,
  2^{2^{\aleph_0}} \,\, = \,\,
  \beth_2 \,\,\, .
\end{equation}
What else has this level of cardinality?
Some say that the set of all possible {\em functions} of real numbers
(i.e., functions that map real numbers onto real numbers) has
cardinality $2^{\mathfrak c}$.
However, it would probably take many more pages to properly specify
which kinds of functions are applicable to this concept, and which
ones are not.
One thing we can say for sure is that the {\em smallest} possible
cardinality that the power-set of the reals can have is $\aleph_2$.

One can keep going up the ladder of higher $n$ for both
$\aleph_n$ and $\beth_n$,
and then take that process to an infinite limit to obtain,
say, $\aleph_{\omega}$ and $\beth_{\omega}$.
Using the set-theory notation discussed earlier, we can define
these more precisely as unions of smaller sets; i.e.,
\begin{equation}
  \aleph_{\omega} \,\, = \,\, \bigcup_{n \in {\mathbb N}} \,
  \aleph_n 
  \,\,\,\,\,\,\,\,\,\,\,\,
  \mbox{and}
  \,\,\,\,\,\,\,\,\,\,\,\,
  \beth_{\omega} \,\, = \,\, \bigcup_{n \in {\mathbb N}} \, \beth_n 
  \,\,\, .
  \label{eq:subomega}
\end{equation}
Note that $\aleph_{\omega}$ is the smallest $\aleph$-number for which
there must exist a countably infinite number ($\omega$) of other
$\aleph$-numbers that are smaller than it.

There's an interesting piece of trivia about $\aleph_{\omega}$.
It has been shown with ZFC that this is the first
uncountable cardinal that {\em cannot} be equal to ${\mathfrak c}$
\citep[see, e.g.,][]{Ea70}.
Unfortunately, this does not put an upper bound on speculations about
the continuum hypothesis, since all it tells us is that
${\mathfrak c} \neq \aleph_{\omega}$.
Thus, it {\em may} be the case that ${\mathfrak c} > \aleph_{\omega}$.

We can finish this section by taking stock of the cardinals that we
know about so far, and we can sort them by whether they are successor
cardinals or limit cardinals (see Section~\ref{sec:countable:lim1}
for this distinction amongst ordinal numbers).
Note, however, that there are now two ways to pin down whether one
cardinal comes ``after'' another.
First, there's the standard successor operation (i.e., going from
$\kappa$ to $\kappa^{+}$).
Second, there's the power-set operation (i.e., going from
$\kappa$ to $2^{\kappa}$).
Both of these operations create a cardinal that is larger than
the original one, and mathematicians often use this distinction to
split up the definition of a limit cardinal into two types:
\begin{itemize}
\item
A {\em weak limit cardinal} has no immediate predecessor, using
the standard direct successor operation.
This is the same meaning of ``limit cardinal'' that extends from
our use of ``limit ordinal'' earlier.
The technical definition is that $\kappa$ is a weak limit cardinal
if and only if, for every other smaller cardinal $\lambda < \kappa$,
it is true that $\kappa > \lambda^{+}$.
If we were able to find a $\lambda$ for which $\kappa = \lambda^{+}$,
then $\lambda$ would be the immediate predecessor of $\kappa$, and
thus $\kappa$ would be a successor cardinal, {\em not} a limit cardinal.
\item
A {\em strong limit cardinal} has no ``power-set predecessor.''
The technical definition is that $\kappa$ is a strong limit cardinal
if and only if, for every other smaller cardinal $\lambda < \kappa$,
it is true that $\kappa > 2^{\lambda}$.
If we were able to find a $\lambda$ for which $\kappa = 2^{\lambda}$,
then $\lambda$ would be the power-set predecessor of $\kappa$, and
thus $\kappa$ would not be a strong limit cardinal.
(It still may be either a weak limit cardinal or a successor cardinal.)
\end{itemize}
Depending on the status of GCH, the size of the ``step'' taken
when doing the power-set is either larger than or equal to the step
taken when doing a successor operation (see Equation~(\ref{eq:bgeqa})).
Thus, all strong limit cardinals must also be weak limit cardinals.

Here is a summary of the status of the cardinals we know about so far:
\begin{itemize}
\item
As before, let us just skip zero because it's a special case.
\item
The finite positive numbers $n = 1, 2, 3, 4, \ldots$ are
successor ordinals, because they all have immediate predecessors.
Each one corresponds to a unique and finite successor cardinal number,
which we use the same numerals to denote.
\item
The first infinite ordinal $\omega$ has no immediate predecessor, so
it is a limit ordinal that corresponds to the (first!)\  limit cardinal
$\aleph_0$.
In fact, this is also a strong limit cardinal because one cannot find
any finite number $n$ for which $2^n$ is equal to $\aleph_0$, either.
\item
$\aleph_1$, $\aleph_2$, $\aleph_3$, and all $\aleph_n$ for finite $n$,
have immediate predecessors, so they are successor cardinals.
\item
$\beth_1$, $\beth_2$, $\beth_3$, and all $\beth_n$ for finite $n$,
are also successor cardinals, but it is not clear where their
predecessors fall in the $\aleph$ hierarchy.
\item
$\aleph_{\omega}$ has no immediate predecessor, so it is a
limit cardinal.
Because we don't know much about how it connects to its neighbors
via the power-set operation, it can only be a weak limit cardinal.
\item
$\beth_{\omega}$ cannot be written as the power-set of any lesser
cardinal, so it is a strong limit cardinal.
\item
${\mathfrak c}$, the cardinality of the continuum, could be a
successor cardinal or a limit cardinal.
If we stay within the axioms of ZFC, we just don't know.
\end{itemize}

\vspace*{0.13in}
\section{Beyond the Alephs?}
\label{sec:large}

We are now entering the domain of the {\em large cardinals}
\citep[see, e.g.,][]{Dr74,Je03,Ka09,Bg13,Ho17,Sr22}.
This is also taking us beyond the usual confines of standard ZFC
set theory, because the actual existence of large cardinals cannot
actually be proven rigorously in ZFC.
This is the realm of active, exploratory research.

\subsection{Initial Concepts}
\label{sec:large:init}

One can conceive of cardinals that break out of the hierarchy
defined above.
For example, we already started to discuss the possible meaning of
$\aleph_{\omega}$.
If we momentarily ignore the distinction between ordinals and
cardinals, we can write this as
\begin{displaymath}
  \aleph_{\aleph_0} \,\,\, ,
\end{displaymath}
and then we can start to think about what {\em its} successor may be.
Because $\omega$ and $\omega + 1$ have the same cardinality,
I don't think it makes a lot of sense to talk about
\begin{displaymath}
  \aleph_{\aleph_0 + 1} \,\,\, .
\end{displaymath}
However, there is one obvious way to write down a cardinal that is most
definitely {\em larger.}
That would be
\begin{equation}
  \aleph_{\aleph_1} \, > \, \aleph_{\aleph_0} \,\,\, .
\end{equation}
We should be clear that $\aleph_{\aleph_1}$ is immensely
larger than $\aleph_{\aleph_0}$, and is not its immediate successor.
Still, this now lets us ``leap'' to still higher cardinalities, with
the sequence
\begin{equation}
  \aleph_{\aleph_0} \,\, , \,\,
  \aleph_{\aleph_1} \,\, , \,\,
  \aleph_{\aleph_2} \,\, , \,\,
  \aleph_{\aleph_3} \,\, , \,\,\, \ldots \,\,\, , \,\,
  \aleph_{\aleph_{\aleph_0}} 
\end{equation}
corresponding to ever-larger cardinal numbers.
Taking that kind of limit is just the first step in an even more
expansive sequence that would let us write down
\begin{equation}
  \aleph_{\aleph_0}
  \,\, , \,\,
  \aleph_{\aleph_{\aleph_0}}
  \,\, , \,\,
  \aleph_{\aleph_{\aleph_{\aleph_0}}}
  \,\, , \,\,
  \ldots
  \,\, , \,\,
  \aleph_{\aleph_{\aleph_{\aleph_{\aleph_{\ddots}}}}} \,\,\, .
\end{equation}
Once we have an infinite chain of subscripts,
adding on another subscript $\aleph$ changes nothing.
We can recall that $\aleph_n$ is defined as the cardinal
preceded by $n$ smaller $\aleph$-numbers.
Thus, our new infinite-chain $\aleph$ number is preceded by a
number of smaller cardinals essentially equal to itself!
We define a cardinal $\kappa$ to be an $\aleph$-fixed point if
\begin{equation}
  \kappa \, = \, \aleph_{\kappa} \,\,\, .
\end{equation}
There is a similar definition for $\beth$-fixed points, which
we can express as
\begin{equation}
  \lambda \, = \, \beth_{\lambda} \,\,\, ,
\end{equation}
and we can extrapolate from Equation~(\ref{eq:bgeqa}) to
conclude that $\kappa \leq \lambda$.

Recall that this is not our first encounter with fixed points.
We saw a bunch of them in Sections~\ref{sec:powertower} and
\ref{sec:countable}.
The last one we saw at the ``end'' of the countable infinities
was $\varepsilon_{\alpha} = \alpha$.
This bears the most resemblance to our current situation, since the
quantity that the condition defines is in the subscript of some
other quantity.

Inspired by \citet{Ru82}, let us assign the general symbol $\theta$ to
this new class of numbers corresponding to $\aleph$-fixed points.
``We're far from the shallow, now,'' but we can dip a toe into into the
deeper waters by contemplating analogies in this chart (slightly
modified from one given by Rucker):
\begin{displaymath}
  \begin{array}{|c@{\hspace*{0.3in}}c@{\hspace*{0.3in}}c@{\hspace*{0.3in}}c@{\hspace*{0.3in}}c|}
  \hline
    0 & 1      & 2        & \cdots & \omega \\
  \hline
    0 & \omega & \aleph_1 & \cdots \, \cdots \, \cdots & \theta \\
  \hline
  \end{array}
\end{displaymath}
Let us discuss going from the 1st column to the 2nd column.
On the top row, this is a transition from nothing (0) to something (1).
On the bottom, it is the transition from finite to infinite;
specifically, to the {\em first} infinite ordinal.
Thus, we can think of the bottom row as an ``infinitely boosted''
analog of the top row.

Going from the 2nd to the 3rd column is a ``successor'' operation:
2 is the successor of 1, and $\aleph_1$ is the successor cardinal of
$\omega$ (i.e., $\aleph_0$).

Then, going from the 3rd column to the final column, we see the
limit of taking an infinite number of successor operations.
The top represents a {\em countably} infinite number of successors,
which gets us to $\omega$.
It may sound like a tautology to say that the only way to reach
$\omega$ is to take $\omega$ successors, but it is useful to think
about it in this way.
When we go to the bottom row, we cannot merely take $\omega$
successors again.
That would get us to $\aleph_{\omega}$.
The ``boosted'' bottom row requires us to consider an
{\em uncountably} infinite number of successors.
For now, I'll state without proof that the inheritor of this crown
is $\theta$.
There are some important caveats to make below, but for now let us just
say that the only way to get to $\theta$ is to take $\theta$ successors.
One cannot reach it by taking any finite (or even any {\em countable})
number of successor operations on an $\aleph$ number.

The top row of Rucker's chart gets us to the first infinite number
$\omega$.
By analogy, the bottom row gets us to the first {\em large cardinal.}

\subsection{Inaccessibility}
\label{sec:large:inaccess}

When I introduced $\aleph_{\omega}$ back in Section~{\ref{sec:beyondC}},
I was surreptitiously introducing a completely new type of cardinal.
Up until the advent of $\aleph_{\omega}$, every new cardinal
that came along was truly ``the first of its kind.''
For example, $\aleph_0$ corresponds to $\omega$, the very {\em first}
countably infinite number.
Then, when we stepped from $\aleph_0$ (countable) to
$\aleph_1$ (uncountable), it was assumed that we took a leap to
the {\em first} uncountably infinite number of this kind.

These ``firsts'' have the property that there is no way to describe
them using only the numbers that came before.
Specifically,
\begin{itemize}
\item
There's no way to juggle a {\em finite} collection of {\em finite}
numbers to completely describe $\omega$.
\item
There's no way to juggle a {\em countably infinite} collection of
{\em countably infinite} numbers to completely describe $\aleph_1$.
\end{itemize}
The only way to describe any such cardinal is by using...\  itself!
Cardinals of the above variety are called {\em regular cardinals.}

Another way to say this is that $\kappa$ is a regular cardinal if
it is impossible to ``reach it'' (using any form of counting)
in less than $\kappa$ steps.

Now, however, there is a new beast in town: $\aleph_{\omega}$.
It is the largest member of the set that contains
\begin{displaymath}
  \aleph_0 \, , \,\,\,
  \aleph_1 \, , \,\,\,
  \aleph_2 \, , \,\,\,
  \aleph_3 \, , \,\,\,
  \aleph_4 \, , \,\,\, \ldots \,\,\,\,
  \mbox{and all $\aleph_n$ (for finite $n$) .}
\end{displaymath}
As we said when discussing Rucker's chart, $\aleph_{\omega}$
involves taking the limit of a series that contains ``only''
$\omega$ members.
In this sense, it {\em is} possible to juggle a countably infinite
collection of items to completely describe
$\aleph_{\omega}$ (which itself is uncountably infinite).
This puts $\aleph_{\omega}$ into a different category than the
other regular cardinals.  We call these {\em singular cardinals.}
All cardinals are either regular or singular.

We can pause to provide a way to define a singular cardinal using
set-theory notation \citep[see, e.g.,][]{Pe10}.
If an infinite cardinal $\kappa$ is singular, then it {\em can} be
written in the following form:
\begin{equation}
  \kappa \,\, = \,\, \bigcup_{n \in \lambda} \mu_{n}
  \,\,\, ,
\end{equation}
where $\lambda < \kappa$ and every member of this union
$\mu_{n} < \kappa$ as well.
Note that Equation~(\ref{eq:subomega}) shows how both
$\aleph_{\omega}$ and $\beth_{\omega}$ satisfy this condition,
since either of them can be written as the union of
a set of members whose cardinality ($\omega$) is {\em less than}
the cardinality of the desired thing itself
($\aleph_{\omega}$ or $\beth_{\omega}$).
On the other hand,
if we ran out of smaller subsets and could only describe some
cardinal $\kappa$ by using a subset of {\em equal} cardinality to
$\kappa$ itself, we would have a regular cardinal.

If readers want to learn more about the distinction between
regular and singular cardinals, they would benefit from reading
about the set-theory property of ``cofinality.''
This property crystallizes the above ideas (i.e., whether or not
we can juggle smaller collections in order to form a complete
description) in even more precise terms.
The cofinality of a cardinal $\kappa$, which we call
$\mbox{cf}(\kappa)$, is the smallest number that provides a
``complete description'' of $\kappa$.
Thus, regular cardinals have $\mbox{cf}(\kappa) = \kappa$,
and singular cardinals have $\mbox{cf}(\kappa) < \kappa$.
It can also be shown (using some parts of ZFC) that
\begin{displaymath}
  \begin{array}{rcl}
  \mbox{If $\delta$ is a successor ordinal (or 0), then} \,\,\, &
  \mbox{cf}(\aleph_{\delta}) \, = \, \aleph_{\delta} \,\,\, , &
  \,\,\, \mbox{and $\aleph_{\delta}$ is regular .} \\
  \mbox{If $\lambda$ is a limit ordinal, then} \,\,\, &
  \mbox{cf}(\aleph_{\lambda}) \, = \, \lambda \, < \, \aleph_{\lambda}
  \,\,\, , &
  \,\,\, \mbox{and $\aleph_{\lambda}$ is singular .}
  \end{array}
\end{displaymath}
One implication of this is that, for any cardinal $\kappa$, its
immediate successor $\kappa^{+}$ must be regular.
This principle covers nearly all cardinals, but it skips over limit
cardinals like $\aleph_{\omega}$, which have no immediate predecessors.
Thus,
\begin{itemize}
\item
If $\kappa$ is a successor cardinal, it is regular.
\item
If $\kappa$ is a limit cardinal, it may be regular or it may be singular.
\end{itemize}
We can provide one additional piece of interesting trivia: recall our
claim from Section~\ref{sec:beyondC} that ${\mathfrak c}$, the
cardinality of the real continuum, cannot be equal to $\aleph_{\omega}$.
We now have the vocabulary to express this in a more general way:
it must be the case that $\mbox{cf}({\mathfrak c}) > \aleph_0$ (i.e.,
that both the continuum and its cofinality must be {\em uncountable}).

How does all this relate to $\theta$, our general term for
$\aleph$-fixed points?
It depends on how we actually define the infinite chain of
subscripts in
\begin{equation}
  \aleph_{\aleph_{\aleph_{\aleph_{\aleph_{\ddots}}}}}
\end{equation}
The most simple (i.e., straightforward) way that we can specify this
chain is to say that it involves adding on a countably infinite ($\omega$)
number of $\aleph$ subscripts.
This, after all, is what we meant by the dots ($\cdots$) going all the
way back to Equation~(\ref{eq:natural}).
Note that this description involves the union of ``only'' $\omega$
members, so it is similar to the example of $\aleph_{\omega}$.
Thus, we must conclude that the {\em smallest} $\aleph$-fixed point
is a singular cardinal.

However, we know that fixed-point conditions allow for multiple solutions.
What about those?
If we counted up a finite number of succeeding solutions to
$\theta = \aleph_{\theta}$, or even a countably
infinite number of them, we are still in the land of singular
cardinals---since their description involves numbers of smaller
cardinality than the cardinal that we're trying to describe.
In order to really define a cardinal that is fundamentally larger than
everything that has come before, we need to take
{\em that number} of succeeding solutions to the fixed-point relation.
\citet{Ta38} proved that the first of these, which we now call $\theta_0$,
is both a limit cardinal and a {\em regular} (not singular) cardinal.
So far, the only other example of this combination has been $\aleph_0$.
Earlier than that, \citet{Hau08} introduced the term {\em inaccessible}
to describe regular limit cardinals, and we add to that by restricting
the definition to all regular limit cardinals greater than $\aleph_0$.

\begin{figure}[!t]
\epsscale{1.17}
\plotone{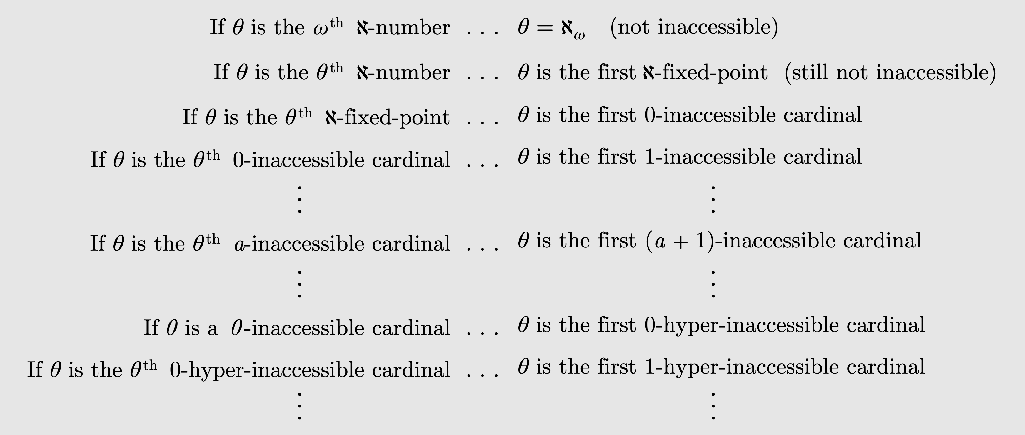}
\caption{Renaming and reframing the first few generations of large cardinals.}
\label{fig:theta_fixed}
\end{figure}

To reiterate: one can only reach the smallest inaccessible cardinal
$\theta_0$ by making an infinitely long list of succeedingly larger
and larger solutions to the $\aleph$-fixed point condition, and not
stopping until one counts $\theta_0$ of those solutions!
Another way of laying out these self-referential conditions is shown
in the first few rows of Figure~\ref{fig:theta_fixed}.

One can also make distinctions amongst the inaccessibles by
defining the weakly-inaccessible cardinals (corresponding to
regular weak-limit cardinals) and the strongly-inaccessible
cardinals (for regular strong-limit cardinals).
Note that all inaccessible cardinals are $\aleph$-fixed points,
but not all $\aleph$-fixed points are inaccessible cardinals.
All strongly-inaccessible cardinals must be both
$\aleph$-fixed points and $\beth$-fixed points,
and if GCH is true then these distinctions disappear.

The reason these quantities are considered to be so very {\em inaccessible}
(i.e., unreachable from lower levels of cardinality) comes from
the combination of the two defining properties:
\begin{itemize}
\item
Because it is regular, it cannot be reached by taking unions of
sets smaller than itself.
\item
Because it is a limit cardinal, it cannot be reached by
``leaping from cardinal to cardinal'' (either via the direct 
successor operation or the power-set operation, depending on
whether it is strong or weak) with any number of leaps smaller
than itself.
\end{itemize}
Essentially, inaccessible cardinals occupy a ``tier'' unto themselves,
in that none of the standard set-theory operations that can be done to
cardinals in the lower (accessible) tiers can produce them.

\subsection{Degrees of Inaccessibility}
\label{sec:large:larger}

Just like there were infinitely many $\varepsilon$ ordinals
(i.e., ordinals that satisfy the fixed-point condition $\omega^a = a$),
it is true that there exist
an infinite number of inaccessible cardinals that
satisfy $\theta = \aleph_{\theta}$, too.
Thus, it may be possible to systematize them as a sequence of
cardinals that we can write as $\theta_n$, for some numerical index
$n$ that sorts them by largeness.
Actually, this enumeration may only cover a subset of
{\em all} of the possible inaccessibles, but as long as it counts
{\em some} of them---and it orders them in cardinality from least to
greatest---it is useful to think about the sequence:
\begin{displaymath}
  \theta_0 \,\, < \,\,
  \theta_1 \,\, < \,\,
  \theta_2 \,\, < \,\,
  \theta_3 \,\, < \,\,\, \cdots \,\,\, < \,\,
  \theta_{\omega} \,\, < \,\,\, \cdots \,\,\, < \,\,
  \theta_{\theta_0} \,\, < \,\,\, \cdots
\end{displaymath}
I have seen it suggested that one way to reach a larger inaccessible
is to start with $\theta_n$ and then take the limit of the sequence
\begin{equation}
  \theta_n \,\, , \,\,\,\,
  \aleph_{\theta_n} \,\, , \,\,\,\,
  \aleph_{\aleph_{\theta_n}} \,\, , \,\,\,\,
  \aleph_{\aleph_{\aleph_{\theta_n}}} \,\, , \,\,\,\, \ldots
\end{equation}
and then call that new cardinal $\theta_{n+1}$.
However, if we recall Section~\ref{sec:powertower},
it may be that there are inaccessible cardinals that are solutions
to the $\aleph$-fixed-point equation but are {\em not}
expressible as infinitely iterated chains of $\aleph$ subscripts.
No matter what method is chosen to find successively larger inaccessibles,
one can eventually reach $\theta_{\theta_0}$, then eventually 
$\theta_{\theta_{\theta_0}}$, and so on.
If that process is taken to its own limit, we obtain yet another
new type of cardinal,
\begin{equation}
  \nu \,\, = \,\, \theta_{\theta_{\theta_{\theta_{\ddots}}}}
\end{equation}
which is now our first example of a {\em 1-inaccessible cardinal.}
This is kind of a bland name, but it implies that what we have been
calling merely inaccessible should have been called {\em 0-inaccessible}
all along.

Whenever we see that infinite chain of subscripts, we should now
be ready to see a new fixed-point relation.
It's true that 1-inaccessible numbers like $\nu$ satisfy such a
relation:
\begin{equation}
  \theta_{\nu} \,\, = \,\, \nu \,\,\, .
\end{equation}
Just like the first 0-inaccessible cardinal $\theta_0$
was the $\theta_0^{\rm th}$ $\aleph$-fixed point,
one can think of $\nu$ as the $\nu^{\rm th}$
0-inaccessible cardinal.
In other words, we define $\nu$ as 1-inaccessible if $\nu$ is
0-inaccessible, and if there are at least $\nu$ 0-inaccessible cardinals 
less than $\nu$.
\citet{Ru82} provided one final way to wrap our heads around this
concept (despite using slightly different terminology).
Recall that 0-inaccessible cardinals were considered ``out of reach''
because of two criteria: they are regular cardinals, and they are
limit cardinals.
Well, a 1-inaccessible cardinal is even {\em more out of reach} because
of three criteria:
\begin{itemize}
\item
Because it is regular, it
cannot be reached by taking unions of sets smaller than itself.
\item
Because it is a limit cardinal, it cannot be reached by leaping from
cardinal to cardinal
(either via $\kappa^{+}$ or $2^{\kappa}$).
\item
Because it is also the limit of an infinite sequence of 0-inaccessibles,
it cannot be reached by ``leaping from inaccessible to inaccessible.''
\end{itemize}
Maybe it is overkill to include all three items in this list, since
the third item requires the first two to be true already.

Of course, once we have 1-inaccessibles as the limit of the infinite
sequence of 0-inaccessibles, we can define the 2-inaccessibles as
the limit of the infinite sequence of 1-inaccessibles.
This is followed by 3-inaccessibles, 4-inaccessibles, and then by
the now-infamous words ``and so on...''\  to
obtain the $(a+1)$-inaccessibles as the limit of the infinite series
of $a$-inaccessibles.
A few of these steps are illustrated in the
middle rows of Figure~\ref{fig:theta_fixed}.

Where does {\em this} railroad take us?
As one steps up the ladders of $a$-inaccessibility, there is an
infinite limit at which adding one more more layer of inaccessibility
changes nothing.
Thus, we define a large cardinal $\kappa$ to be
{\em hyper-inaccessible} if it is $\kappa$-inaccessible.
(Maybe another way to refer to it is ``itself-inaccessible!'')
This is a new fixed-point that is a bit more difficult to express
as an equation.
Let us attempt to do that by defining some new terminology.
The 0-inaccessibles can be called $\theta_{n}(0)$,
then the 1-inaccessibles are $\theta_{n}(1)$, and the general
$a$-inaccessibles are $\theta_{n}(a)$. 
Then, we can see that
\begin{equation}
  \mbox{if}
  \,\,\,\,\,\,\,
  \theta_{\kappa}(a) \, = \, \kappa  \,\, ,
  \,\,\,\,\,
  \mbox{then $\kappa$ is $(a+1)$-inaccessible, and we write it as}
  \,\,\,\,\,
  \theta_{n}(a+1) \,\, .
  \label{eq:a_inacc_fixed}
\end{equation}
Note that the new notation $\theta_{n}(a+1)$ has an undetermined
subscript $n$, which indexes the infinite number of inaccessibles
at this new tier.
Now, after repeating the above fixed-point procedure an infinite
number of times, we reach the hyper-inaccessibles.
For them,
\begin{equation}
  \mbox{if}
  \,\,\,\,\,\,\,
  \theta_{n}(\kappa) \, = \, \kappa  \,\, ,
  \,\,\,\,\,
  \mbox{then $\kappa$ is hyper-inaccessible, and we write it as}
  \,\,\,\,\,
  \theta_{n}(0,1) \,\, .
  \label{eq:hyper_inacc_fixed}
\end{equation}
Please don't move on before taking note of the subtle indexing shifts
that happen as one goes from Equation~(\ref{eq:a_inacc_fixed}) to
Equation~(\ref{eq:hyper_inacc_fixed}).
In the mathematical literature, there are quite a few alternate
definitions of the different levels of inaccessibility and
hyper-inaccessibility, and I cannot guarantee that the above
corresponds to the most commonly found versions.

I hate to keep rushing down the rabbit-hole without much of a pause
for breath, but it is possible to generalize even more.
What we called hyper-inaccessible above can be called
0-hyper-inaccessible.
Then, we can define $\kappa$ as a 1-hyper-inaccessible cardinal if
$\kappa$ is already 0-hyper-inaccessible and if there are at least
$\kappa$ 0-hyper-inaccessible cardinals less than $\kappa$.
The same limit-taking gives us 2-hyper-inaccessible cardinals,
3-hyper-inaccessible cardinals, and so on.
We can unify all of these tiers using the notation defined above, as
\begin{equation}
  \begin{array}{l}
    \theta_{n} (i,0) \,\,\,
    \mbox{for the $i$-inaccessibles,
    which we used to call $\theta_n(i)$ .} \\
    \theta_{n} (i,1) \,\,\,
    \mbox{for the $i$-hyper-inaccessibles .} 
  \end{array}
\end{equation}
Then, as we climb the ladder of $i$-hyper-inaccessibles, we can take
the limit once again to define a large cardinal $\kappa$ to be
{\em hyper-hyper-inaccessible} (or hyper$^2$-inaccessible) if it
is $\kappa$-hyper-inaccessible.
Of course, we can then go on to define hyper$^3$-inaccessible,
hyper$^4$-inaccessible, and so on, and the notation generalizes to
\begin{equation}
  \begin{array}{rl}
    \theta_{n} (i,h) \,\,\,\, &
    \mbox{for the $i$-hyper$^{h}$-inaccessibles .} 
  \end{array}
\end{equation}
Note that the 0-hyper$^0$-inaccessibles are just the original
``plain'' inaccessibles.

\citet{Cm15} outlined how we can keep climbing up the ladder of
inaccessible cardinals.
We saw that when the $i$-index is taken to an infinite fixed-point limit
(i.e., when we defined a new cardinal $\kappa$ as being
$\kappa$-inaccessible), we were uplifted from the inaccessibles to
the hyper-inaccessibles.
The $h$-index can also be taken to an infinite fixed-point limit,
and that means we can define a large cardinal $\kappa$ that is
hyper$^{\kappa}$-inaccessible.
\citet{Cm15} called this new tier the {\em richly inaccessible cardinals,}
and our notation blossoms forth with a third index in parentheses; i.e.,
we can use
\begin{equation}
  \begin{array}{rl}
    \theta_{n} (i,h,r) \,\,\,\, &
    \mbox{for the $i$-hyper$^{h}$-richly$^{r}$-inaccessibles ,} 
  \end{array}
\end{equation}
where the smallest richly inaccessible cardinal would be called
$\theta_{0}(0,0,1)$.
\citet{Cm15} found useful results by rocketing up at least five additional
tiers and defining the {\em utterly, deeply, truly, eternally,} and
{\em vastly inaccessible} cardinals, which we can go on to specify as
\begin{equation}
  \theta_{n} (i,h,r,u,d,t,e,v) \,\,\, .
\end{equation}
Of course, we didn't have to stop there, but it probably makes sense
to look at things in a slightly different way.

\subsection{Top of the Heap?}
\label{sec:large:mahlo}

\citet{Ma11} defined a class of cardinals that sits at the end of the
above conceptual chain of inaccessibles.
Such cardinals, for which some use the notation $\rho$, are now
called {\em Mahlo cardinals,} and they must be 
$\rho$-inaccessible, hyper$^{\rho}$-inaccessible,
richly$^{\rho}$-inaccessible, utterly$^{\rho}$-inaccessible,
deeply$^{\rho}$-inaccessible, truly$^{\rho}$-inaccessible,
eternally$^{\rho}$-inaccessible, vastly$^{\rho}$-inaccessible, 
and so on for {\em all} subsequent tiers.
In fact, for any $\rho$ that is a Mahlo cardinal, there must also exist
$\rho$ smaller cardinals (call them $\kappa$, with all $\kappa < \rho$)
that are $\rho$-inaccessible, hyper$^{\rho}$-inaccessible,
richly$^{\rho}$-inaccessible, and so on, with all of those
smaller $\kappa$'s not being Mahlo.
Thus, although a Mahlo cardinal must satisfy all of the above
inaccessibility properties, not all cardinals that satisfy these
properties are Mahlo cardinals.

\citet{Ru82} extended the two-row table that we introduced
in Section~\ref{sec:large:init} in order to help us wrap our heads
around how large Mahlo cardinals really are:
\begin{displaymath}
  \begin{array}{|c@{\hspace*{0.3in}}c@{\hspace*{0.3in}}c@{\hspace*{0.3in}}c@{\hspace*{0.3in}}c@{\hspace*{0.3in}}c@{\hspace*{0.3in}}c|}
  \hline
    0 & 1      & 2        & \cdots & \omega & \longrightarrow & \aleph_1 \\
  \hline
    0 & \omega & \aleph_1 & \cdots \, \cdots \, \cdots & \theta &
    \longrightarrow \, \longrightarrow \, \longrightarrow & \rho \\
  \hline
  \end{array}
\end{displaymath}
Recall, back in Section~\ref{sec:countable}, that we experienced a
somewhat exhausting process of building up more and more complicated
(countable) ordinals that started from $\omega$.
This process never really got us to the next (uncountable) number
$\aleph_1$ until we conceptualized things in a fundamentally new way.
Well, here we built up more and more complicated degrees of
inaccessibility, hyper-inaccessibility, and whatever-inaccessibility,
starting with our first discussions of $\theta$, but it never quite
got us to that next fundamental stage, which we now call $\rho$.

\citet{Ru82} used one additional way to describe how $\rho$ is so much
larger than what came before.
This analogy may come the closest to the precise set-theory definition
of Mahlo cardinals.
(Note that I haven't given this definition, because I don't want to
parrot concepts that I still do not really understand!)
Despite being inaccessible (and thus regular), $\rho$ really cannot be
thought of as ``the first of its kind'' in any concrete way.
Try to envision any imaginable fixed-point property $P$ that an
infinite cardinal might obey.
We've seen quite a few examples of these so far, and there must be an
infinite number of additional examples.
For any $P$ that you can envision,
if a Mahlo cardinal $\rho$ obeys $P$, then there always must exist
a smaller cardinal $\kappa < \rho$ that also obeys $P$.
This makes it immensely more difficult to ``reach'' $\rho$ from below.

Despite how far we have climbed, it is still a possibility
(if GCH is not valid) that the cardinality of the real continuum
$\mathfrak c = 2^{\aleph_0}$ may in fact be sitting at this height;
i.e., it may be of a comparable magnitude as $\rho$.

\section{Further Up the Large Cardinal Ladder}
\label{sec:morelarge}

\subsection{Axioms and Consistency Strength}
\label{sec:morelarge:axioms}

Recall, back in Section~\ref{sec:finite}, we quoted
\citet{Rayo19}, who said that
``Our quest to find larger and larger numbers has now morphed into
a quest to find more and more powerful languages!''
From this point on, the study of large cardinals definitely echoes
this quest, but the languages in question are the highly abstract
logical systems of set theory.

It has been known since the 1960s \citep[see, e.g.,][]{Ko14}
that climbing the ladder to
define ever-larger cardinals requires adding {\em more and more new axioms}
to the basics of ZFC set theory.
These new axioms appear to sort themselves into a strict linear order
in terms of their ``consistency strength.''
Consider a given axiom $A_1$ that produces a theory that is both
self-consistent {\em and} can be used to prove that the theory given
by axiom $A_2$ is also self-consistent.
However, axiom $A_2$ is not able to prove that the theory given by
axiom $A_1$ is self-consistent.
This means that axiom $A_1$ is stronger than axiom $A_2$,
and also its corresponding large cardinals occur (usually) higher up
on the ladder.

\begin{figure}[!p]
\epsscale{0.67}
\plotone{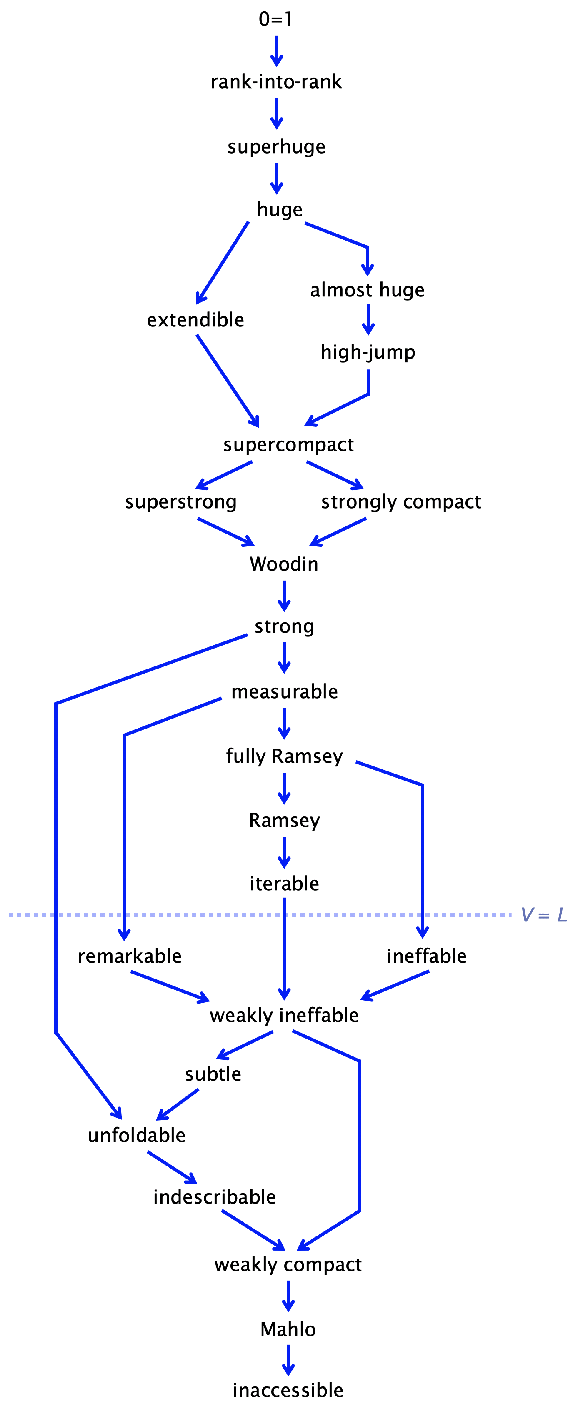}
\caption{The author's idiosyncratic diagram of the large cardinal
hierarchy, with the strongest consistency strengths at the top.
Blue arrows point from theories that directly imply other theories.
\label{fig:largechart}}
\end{figure}

Because the ``higher'' theories are bound by more and more arbitrarily
chosen axioms, they are restricted to ever-narrower domains of set theory.
However, over the years mathematicians have discovered new interconnections
between the ``lower'' levels, which were not evident until they climbed
higher up and expanded their view.
The hope is that these interconnections will eventually knit themselves
together into a denser cluster of roots that will provide more and more
evidence (circumstantial though it may be) for the plausibility of the
whole edifice.
Quoting \citet{G47},
\begin{quotation}
\noindent
There might exist axioms so abundant in their verifiable consequences,
shedding so much light upon a whole discipline, and furnishing such
powerful methods for solving given problems (and even solving them, as
far as that is possible, in a constructivistic way) that quite
irrespective of their intrinsic necessity they would have to be assumed
at least in the same sense as any well-established physical theory.
\end{quotation}
Figure~\ref{fig:largechart} shows my attempt to synthesize several
diagrams found in books \citep[e.g.,][]{Ka09},
papers \citep[e.g.,][]{NW18}, and online
that illustrate the hierarchy of the large cardinals
and their associated axioms.
I cannot vouch for its exactness or completeness.
Note that the inaccessible and Mahlo cardinals discussed above are
at the very bottom.  There's much more beanstalk to climb.

\subsection{Formal Languages and Hierarchies}
\label{sec:morelarge:formal}

Defining the next level of large cardinals requires discussing the
formal language of set theory that is used to constrain their properties.
One can reach up to larger cardinals by using more sophisticated
languages, and we will provide only the roughest outline of what is
meant by ``more sophisticated.''
For more comprehensive introductions to how these kinds of formula
descriptions work, see texts such as \citet{Su13} and \citet{Nw22}.

\citet{Lv65} first outlined a hierarchy of formula classification
that was then extended in an attempt to describe the indescribable.
We can start with a few simple examples of statements in this formal
language, such as:
\begin{equation}
  \forall n \, \exists m \,\, ( m = n + 1 )  \,\,\, ,
\end{equation}
which means: for all values of a natural number $n$, there exists
a natural number $m$ that is equal to $n+1$.
This is a true statement, since it could be stated with more brevity
that ``every natural number has a successor.''

However, note how the meaning changes when the order of the
unbounded universal quantifier ($\forall$, ``for all'') is swapped
with the unbounded existential quantifier ($\exists$, ``there exists''):
\begin{equation}
  \exists m \, \forall n \,\, ( m = n + 1 )  \,\,\, .
\end{equation}
This statement means: there exists a natural number $m$, such
that for every natural number $n$, it is true that $m = n+1$.
This is false, since it is claiming that there is some number $m$
that is the immediate successor of {\em every} natural number!

Another example helps us come up with another way to describe the
first infinite ordinal $\omega$:
\begin{equation}
  \forall x \, \exists y \,\, ( x \in y )  \,\,\, ,
  \label{eq:xyomega}
\end{equation}
which means that when we consider all possible sets corresponding to the
natural numbers ($x$), there exists a set ($y$, which we call $\omega$)
that all of those natural numbers are contained in.
Is it true?  Well, it is if you accept the existence of infinite sets.

In this case, swapping the order of the quantifiers would say essentially
the same thing, but with a different subject/object emphasis:
\begin{equation}
  \exists y \, \forall x \,\, ( x \in y ) \,\,\, ,
\end{equation}
which means that there exists some set ($y$) that contains every
natural number ($x$).
However, starting again with Equation~(\ref{eq:xyomega}) and
swapping the orders of $x$ and $y$ in both places produces a
false statement.  For example, consider
\begin{equation}
  \exists y \, \forall x \,\, ( y \in x )  \,\,\, .
\end{equation}
This means that there exists something ($y$) that is a member of
every possible set ($x$), which is not true.
In some versions of set theory, the ``opposite'' of
Equation~(\ref{eq:xyomega}) is used to define the concept of the empty set:
\begin{equation}
  \forall x \, \exists y \,\, ( x \not\in y )  \,\,\, ,
  \label{eq:xyempty}
\end{equation}
where the symbol $\not\in$ denotes ``not a member of.''
The above statement says there must exist a set
($y$, which we call $\emptyset$) that contains {\em none} of
the universe of other possible sets ($x$).

We want to be able to classify statements in formal language on
the basis of their complexity; i.e., how well they can
convey increasingly intricate and multilayered concepts.
There are two general ways that a statement can be made more complex:
(1) its chain of alternating, nested quantifiers can be made longer,
or (2) each of its variables can be made to represent more and more
complex ideas.
Let us examine both options.

\subsubsection{Language Hierarchy: Quantifier Alternation}
\label{sec:morelarge:formal:1}

\citet{Lv65} noted that longer chains of nested $\exists$
and $\forall$ quantifiers provide more powerful ways of describing
complex ideas.
Thus, it is possible to classify a statement by counting the number
of {\em alternations} between quantifiers that occur:
\begin{itemize}
\item
The simplest statements contain no unbounded quantifiers, and they
are classified as being of type $\Sigma_0$ or $\Pi_0$.
(These two symbols are equivalent to one another, and are sometimes
called by the single name $\Delta_0$.)
\item
$\Sigma_1$ statements consist of any number of $\exists x$ quantifiers,
followed by a $\Pi_0$ statement.
\item
$\Pi_1$ statements consist of any number of $\forall x$ quantifiers,
followed by a $\Sigma_0$ statement.
\item
In general, $\Sigma_{n+1}$ statements consist of any number of
$\exists x$ quantifiers, followed by a $\Pi_n$ statement.
\item
In general, $\Pi_{n+1}$ statements consist of any number of
$\forall x$ quantifiers, followed by a $\Sigma_n$ statement.
\end{itemize}
Above, we saw how the exact ordering of $\exists$ and $\forall$ quantifiers
makes a difference in the meaning of a statement, so it makes sense that
it matters when making these classifications, too.

Some examples of statements that correspond to different levels of the
L\'{e}vy hierarchy, which we provide only in words
(see, e.g., \citeauthor{Mk21} [\citeyear{Mk21}] for details), are:
\begin{itemize}
\item
``$x$ is an ordinal'' can be written as a
$\Sigma_0$ or a $\Pi_0$ statement.
\item
``$x$ is countable'' can be written as a
$\Sigma_1$ statement.
\item
``$x$ is a cardinal'' can be written as a
$\Pi_1$ statement.
\item
``There exists an inaccessible cardinal'' can be written as a
$\Sigma_2$ statement.
\item
The generalized continuum hypothesis (GCH) can be written as a
$\Pi_2$ statement.
\end{itemize}

\subsubsection{Language Hierarchy: Logical Order}
\label{sec:morelarge:formal:2}

There is another independent hierarchy that describes how to
create increasingly complex statements: the ``order''
of predicate logic that is being used to define the variables
$x$, $y$, and so on.
Our example statements above were all taken from
{\em first-order logic,} in which the variables can only be
straightforward ``individuals'' in the adopted domain
(in this case, either numbers or sets).
The most basic form of ZFC discussed earlier is a first-order logic,
but there are some versions that allow the variables to only be
numbers and not sets.

However, there also exist {\em second-order} logical systems,
in which the variables can be either individuals or relations
between those individuals.
In second-order logic, one can use more abstract quantifiers such
as $\forall f(x)$ (for all functions of $x$) or
$\exists C$ (there exists a class of objects, or a set of sets).
Specifying either $\forall$ for $\exists$ on the power-set
${\mathscr P}(S)$ of some unknown group of sets $S$ is also
straightforward to do in second-order logic.

Moving on even higher-order logic systems is less common to see, but
they exist.
Third-order logic, for example, allows one to define variables for
functions of functions, properties of properties, or classes of classes.
Although there are cases in which equivalent statements can be rephrased
in clever ways that allow one to ``translate up or down an order''
\citep[see, e.g.,][]{Va21},
there are {\em some} things that just cannot be stated in a lower-order
logic that become easy to discuss in a higher-order logic.

The basic L\'{e}vy hierarchy, which classifies statements on the
basis of the number of quantifier alternations, can be extended
to also express the logical order of those statements.
This now involves adding a superscript $m$, which corresponds to the
statement being in $(m+1)^{\rm th}$-order logic.
Thus, the general classifications are $\Sigma^m_n$ and $\Pi^m_n$, with
\begin{itemize}
\item
$\Sigma^0_n$ or $\Pi^0_n$ corresponding to statements in first-order logic.
This defines a system called the {\em arithmetical hierarchy.}
\item
$\Sigma^1_n$ or $\Pi^1_n$ corresponding to statements in second-order logic.
This defines a system called the {\em analytical hierarchy.}
\end{itemize}
A statement that can be classified as both
$\Sigma^m_n$ and $\Pi^m_n$ is sometimes called $\Delta^m_n$.
In general, going up one step in order number $m$ is interpreted as
making a ``bigger'' leap in complexity than just going up one step in
L\'{e}vy subscript $n$.
In fact, one way to think about going up to the next order is to
take the limit of the subscript to (countable) infinity:
\begin{displaymath}
  \Sigma^m_{\omega} \,\, = \,\,
  \Pi^m_{\omega} \,\, = \,\,
  \Delta^m_{\omega} \,\, = \,\,
  \Sigma^{m+1}_{0} \,\, = \,\,
  \Pi^{m+1}_{0} \,\, = \,\,
  \Delta^{m+1}_{0} \,\,\, ,
\end{displaymath}
although this may not be rigorously true in all situations.

\subsection{Indescribable Cardinals}
\label{sec:morelarge:indesc}

The introduction of $\Sigma^m_n$ and $\Pi^m_n$ notation allows us to
make some definitive statements about different kinds of infinite
cardinals that can be described (or not!) using a given level
of language in this hierarchy.
\citet{HS61} first made this connection by defining 
{\em $Q^m_n$-indescribable cardinals} as cardinals that cannot be
described by any statement made at level $Q^m_n$ of the hierarchy,
where $Q$ could be $\Sigma$ or $\Pi$ or $\Delta$.
We will see that this language can be used to classify a number of
different kinds of cardinals at different levels of the tree
shown in Figure~\ref{fig:largechart}.

Note that if a cardinal is, say, $Q^m_n$-describable, then it must
be describable for all higher values of both $m$ and $n$, too.
Going to either higher $m$ or $n$ opens up one's universe of
expression, but it also always encompasses what came before.
However, if a cardinal is $Q^m_n$-indescribable, then two things
must be true: (1) it is also indescribable at all lower values of
$m$ and $n$, and (2) it {\em may}
become describable at a higher $m$ or $n$ (or it may not).
The larger the indices get for a $Q^m_n$-indescribable cardinal,
the more difficult it becomes to describe it using finite language,
and thus its cardinality tends to be larger.

Here is a summary, expressed in words and not in formal language
\citep[for which, see][]{Dr74,Je03,Ka09},
of how we can begin to wrap our heads around these different levels:
\begin{itemize}
\item
The ``smallest'' class of infinite cardinals,
like $\aleph_0$, $\aleph_1$, and $\aleph_{\omega}$, were discussed in
Sections~\ref{sec:uncountable} and \ref{sec:beyondC}.
These can all be called {\em accessible}
because they are all $\Pi^0_n$-describable, for some $n$.
Our first-order definition of $\omega$, given in
Equation~(\ref{eq:xyomega}), is an example of this describability.
\item
Next, we make the leap to inaccessible cardinals (Section~\ref{sec:large})
and we can state that every inaccessible cardinal is
$\Pi^0_n$-indescribable, for all $n > 1$.
As we mentioned earlier, large cardinals like the inaccessibles
cannot be proven to exist in standard ZFC set theory.
Thus, it makes sense that these quantities are not completely
describable using a first-order logic system like ZFC
(i.e., $m=0$ in the superscript).
It has also been shown that inaccessibility is equivalent
to being $\Sigma^1_1$-indescribable.
\item
Strongly-inaccessible cardinals are $\Pi^1_0$-indescribable.
If GCH is not true, then strongly-inaccessible cardinals are larger
than weakly-inaccessible cardinals.
\item
\citet{ET61} introduced the {\em weakly compact cardinals,} which
sit directly above the Mahlo cardinals in Figure~\ref{fig:largechart}.
These correspond precisely to being $\Pi^{1}_{1}$-indescribable,
and they can be considered the limit of cardinals that are Mahlo,
$n$-Mahlo, hyper-Mahlo, richly-Mahlo, and so on (although we never
actually defined these higher levels of Mahlo-ness).
\item
\citet{JK69} introduced the {\em ineffable cardinals,} 
which generally correspond to being $\Pi^1_2$-indescribable.
\item
\citet{Ul30} introduced the {\em measurable cardinals,} which we now
know sit higher up in Figure~\ref{fig:largechart} than
many of the ones discussed above.
If $\kappa$ is measurable, then it can be considered the
$\kappa^{\rm th}$ weakly compact cardinal.
All measurable cardinals are $\Pi^2_1$-indescribable.
Note how even third-order ($m=2$) logic isn't up to the task of
describing cardinals at this level!
\item
{\em Totally indescribable cardinals} are
$\Pi^m_n$-indescribable for every finite natural number $m$ and $n$.
They are also $\Sigma^m_n$-indescribable and $\Delta^m_n$-indescribable.
\end{itemize}
Note that the above list skips several steps in Figure~\ref{fig:largechart}.
However, it is often the case that if cardinal $A$ is directly above
cardinal $B$ in the ladder of consistency strength, one can say
that all cardinals of type $A$ are also of type $B$.

\subsection{Another Visit with the Continuum}
\label{sec:morelarge:contin}

As we continue to climb past landmarks on the mountain of
cardinality, we can still wonder about where $\mathfrak c$,
the cardinality of the real numbers, may be.
Back in Section~\ref{sec:uncountable:C} we outlined some
plausible-sounding arguments in favor of ${\mathfrak c} = \aleph_1$
(i.e., the continuum hypothesis).
However, now let us discuss a few additional concepts that {\em may}
support the idea of $\mathfrak c$ sitting in the neighborhood of
the indescribables that we have just, um, described.

\begin{figure}[!t]
\epsscale{1.13}
\plotone{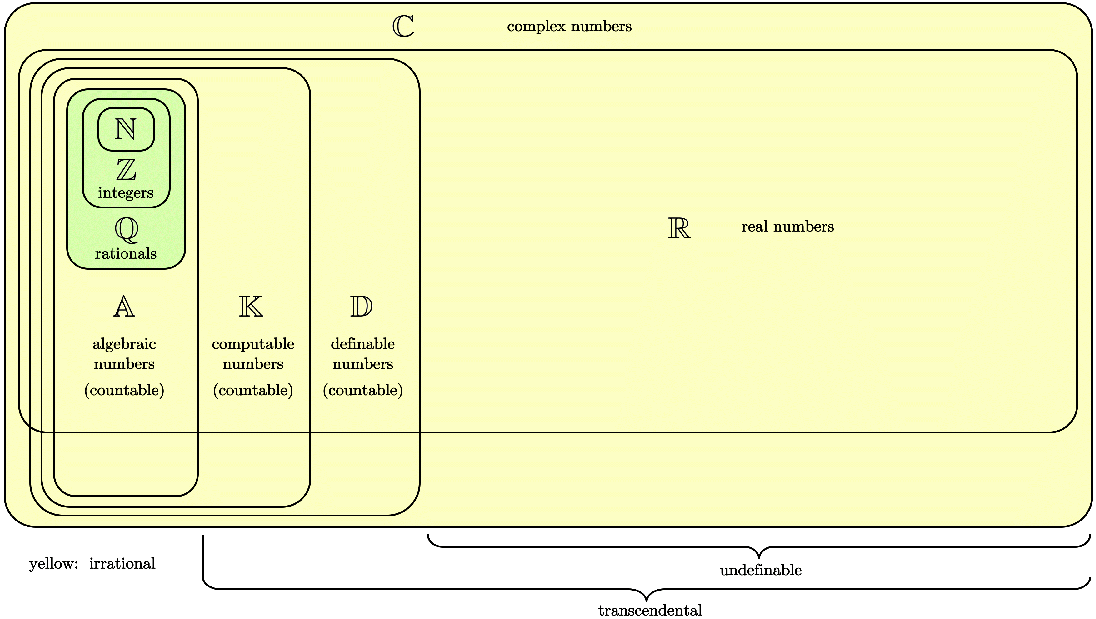}
\caption{Groupings and nested subsets of different types of numbers.
Note that even the smallest set (the natural numbers $\mathbb N$) contains
a countable infinity of members.}
\label{fig:venn_numbers}
\end{figure}

Take a look at the diagram shown in Figure~\ref{fig:venn_numbers}.
We started this whole journey with $\mathbb N$, the set of
all natural numbers.
That was just a subset of the positive and negative integers
(which use the set symbol $\mathbb Z$),
and the integers are a subset of the rationals
(which use the symbol $\mathbb Q$).
This region of {\em countably infinite} sets is outlined in green.
Everything outside that region (shown in yellow) corresponds to
the set of irrational numbers that are either
real (with symbol $\mathbb R$) or complex (with symbol $\mathbb C$),
with both of those full sets having uncountable cardinality $\mathfrak c$.
There are a few interesting sets of numbers in between these extremes:
\begin{itemize}
\item
The {\em algebraic numbers} (for which we use the symbol $\mathbb A$)
consist of the full set of roots of polynomials with integer
coefficients.
Thus, all rational numbers are algebraic, but there are many
irrationals that are also algebraic.
Note that $\sqrt{2}$ is a root of the polynomial $x^2 - 2$, so
this famous irrational number is a member of $\mathbb A$.
Some of these roots are imaginary and complex, so note that
the set $\mathbb A$ crosses over into $\mathbb C$.
We can recall from Section~\ref{sec:countable:lim2}
that the number of polynomials with integer coefficients is
essentially $\omega^{\omega}$, and this corresponds to
cardinality $\aleph_0$.
Thus, despite the algebraic numbers sitting in the yellow realm
that includes the full set of uncountable irrationals,
the set $\mathbb A$ by itself is countable.
The algebraic numbers can fit into Hilbert's Hotel!
\item
The next larger subset shown in Figure~\ref{fig:venn_numbers}
corresponds to the {\em computable numbers,} for which I've seen
people use the symbol $\mathbb K$.
These are numbers that can be computed to arbitrary precision
by use of a finite algorithm.
All algebraic numbers are computable, but not all computable
numbers are algebraic.
Most famously, transcendental numbers such as $\pi$ and $e$
are computable, as one can show by programming a computer to
spit out as many decimal digits of these numbers as one wants.
\citet{Tu36} showed that the set $\mathbb K$ is countable.
\item
Lastly, there is the set $\mathbb D$ of {\em definable numbers.}
These are somewhat similar to the computable numbers, but
here the specification involves being able to describe the
number precisely using formal mathematical language.
If one can write down a formal expression that defines a
number, then it is definable.
One can also show that all computable numbers are definable,
but not all definable numbers are computable.
It is not easy to come up with examples of definable but uncomputable
numbers, but see \citet{Si05} for the fascinatingly weird case of
``Specker sequences.''
Although the possible number of formal expressions one can write
is infinite, it is nevertheless a {\em countable} infinity.
Thus, the set $\mathbb D$ is countable, with cardinality $\aleph_0$,
and its members can fit into Hilbert's Hotel, too.
\end{itemize}
Despite the sets $\mathbb A$, $\mathbb K$, and
$\mathbb D$ covering larger and larger parts of the uncountable
set of irrational numbers, each of those sets by itself is
countable.
There are arguments given in Appendix~\ref{sec:appA} for why the
set of all rationals $\mathbb Q$ takes up negligibly small space
in the full set of reals $\mathbb R$.
However, we can create a one-to-one correspondence between the
members of any countably infinite set.
Thus, we can replace $\mathbb Q$ in these arguments by, say,
$\mathbb N$, $\mathbb Z$, $\mathbb A$, $\mathbb K$, or $\mathbb D$.
This means that the huge set of definable numbers, which are all
numbers that can be described by the formal language of set theory,
is still a {\em negligibly tiny part of the real continuum.}
In this sense, the full diversity of numbers contained in $\mathbb R$
is indescribable using formal set theory!

Does the above set of arguments imply that we can jump to the
conclusion that $\mathfrak c$ is an indescribable cardinal?
Not really.
Although the total number of possible formulas that can be written
down is ``only'' of cardinality $\aleph_0$, we also know that
{\em some} of those formulas can describe infinite cardinalities
larger than $\aleph_0$.
In any case, this line of thinking shows how we may start to
come to terms with the CH being false---and with the possibility that 
$\mathfrak c$ {\em could}
be a much, much larger cardinal than $\aleph_1$.

\section{How Much Higher Can We Climb?}
\label{sec:evenhigher}

\vspace*{0.07in}
\subsection{The Constructible Universe}
\label{sec:evenhigher:VL}

\vspace*{0.07in}
We can keep ascending the ladder of large cardinals shown in
Figure~\ref{fig:largechart}, but it becomes increasingly difficult
to talk about them using plain language.
Earlier, I quoted \citet{Ru82}, who discussed that \citet{G40}
showed how Cantor's continuum hypothesis (CH) is consistent with ZFC.
This did not prove that ${\mathfrak c} = \aleph_1$, but it did show
that one can never prove that ${\mathfrak c} \neq \aleph_1$.
However, this work had additional consequences for whether or not
large cardinals can exist past a certain point, so we need to
talk more about it.

Up until this point, nearly everything we have been discussing has
applied to what is called the {\em von~Neumann universe} of set theory.
The term ``universe'' applies to the collection of all sets that
satisfy a given property.
The von~Neumann universe, which is called $V$, starts with the
sets corresponding to the natural numbers, which we illustrated in
Equation~(\ref{eq:vonNeumann}).
Note how each successor set encompasses all subsets of the preceding set.
$V$ is the union of all of these sets, as well as all other sets
corresponding to infinite successor ordinals (which have immediate
predecessors) and infinite limit ordinals (which don't have immediate
predecessors).
One thing that we can say about $V$ is that there are {\em a lot} of
potential ways to reach some arbitrarily high level of cardinality.
It's also amenable to adding lots of new, exploratory axioms, as we
have already discussed.
Thus, $V$ is kind of a wide-open ``wild west'' of set theory.

\citet{G40} decided to try a minimalist thought-experiment.
Maybe, he thought, it is possible to wrap our heads around things
more easily by defining a parallel universe of set theory that is
{\em more restricted} in how one is allowed climb to the higher levels.
In other words, if we limit ourselves to stricter definitions (i.e.,
of how the higher levels are constructed on the foundations of the
lower levels), then the results may be less dependent on the whims of
theorists who like to play around with axioms.
Thus, G\"{o}del defined a {\em constructible universe} and used
the symbol $L$ for it (for ``Law!'').
The strictness comes into play in the definition of what we mean by
a successor.
Our existing definition of a successor essentially tops out at the
power-set; i.e., we're guaranteed to leap up to a higher level of
cardinality by taking the set of all possible subsets
of what came before.
G\"{o}del replaced ${\mathscr P}(x)$ by a more rigidly definable
operator ${\mathscr D}(x)$, which is the same as the power-set for
finite sets, but is limited to what can be specified by
$\Delta^0_0$-type formulas (using only first-order logic)
for infinite sets.
Strictly speaking, ${\mathscr D}(x)$ specifies a subset of
${\mathscr P}(x)$, so that means the entire universe $L$ must be
a subset of $V$.

The question then becomes: How far can we go by assuming that $V=L$?
It would be nice if the ${\mathscr D}(x)$ successor operation was
all that we really needed.
Also, it would be nice if we were able to extract all of the useful
results of infinite set theory by staying inside the safe
confines of $L$; i.e., a universe that is so
much easier to work in, because of its restrictedness.
\citet{G40} postulated an {\em Axiom of Constructibility} that
postulated that $V=L$, and here are some consequences of that assumption:
\begin{itemize}
\item
When working in ``plain'' Zermelo–Fraenkel theory (ZF), it is not
possible to disprove the Axiom of Choice.
Thus, if ZF is self-consistent, then so is ZFC.
\item
Both the CH and GCH are true, so $\mathfrak c = \aleph_1$ after all
(in the $L$ universe).
\item
Thus, all weak-limit cardinals are also strong-limit cardinals.
\item
However, in $L$, there exists a boundary in consistency strength
(see the horizontal line labeled $V=L$ in Figure~\ref{fig:largechart})
above which large cardinals {\em cannot exist.}
Inaccessible, Mahlo, and weakly compact cardinals sit below this
boundary, but measurable cardinals sit above it.
\end{itemize}
Thus, $L$ has some benefits, but it also has limitations.
In order to go above the $V=L$ line, one has to posit new axioms
that allow for the existence of ``non-constructible'' sets; i.e.,
sets that wouldn't be able to exist in $L$.
For this reason, quite a few mathematicians believe that the
``universe'' is quite a bit larger than $L$.

\subsection{Reflection from the Ultimate Heights}
\label{sec:evenhigher:RP}

To climb higher, I think we need
to discuss one additional concept: {\em the Reflection Principle.}
To begin discussing it, though, we need to begin wrapping our heads
around the idea of Absolute Infinity...\  i.e.,
the cardinality associated with the utter totality of the
von~Neumann universe $V$.
The Absolute is the highest degree of infinity, which by
definition includes all other possible ways of conceptualizing it.
Cantor contemplated this idea with earnestly religious reverence.
Using one interpretation of Cantor's writings, one can use the
symbol $\Omega$ to refer to the Absolute as the union of all ordinals,
and \textcjheb{t} (the final Hebrew letter, ``tav'') as the union of
all cardinals.

The Reflection Principle can be stated as the claim that, for any
statement that can be made about the von~Neumann universe $V$, it is
possible to find a subset of $V$ for which that statement is true.
\citet{Ru82} provided the following explanation:
\begin{quotation}
\noindent
The motivation behind the Reflection Principle is that the Absolute should
be totally inconceivable. Now, if there is some conceivable property $P$
such that the Absolute is the only thing having property $P$, then I can
conceive of the Absolute as ``the only thing with property $P$.''
The Reflection Principle prevents this from happening by asserting that
whenever I conceive of some very powerful property $P$, then the first
thing I come up with that satisfies $P$ will {\em not} be the Absolute,
but will instead be some smallish rational thought that just happens to
reflect the facet of the Absolute that is expressed by saying it
has property $P$.
\end{quotation}
According to its proponents, if we deny the Reflection Principle, we
are claiming that the Absolute can be ``pinned down by finiteness.''
If we could find some finite thing to say that defines the Absolute
(by applying only to it, and not to any smaller subsets of it), then
what we {\em thought} was the Absolute couldn't have really been
the Absolute!
\citet{We19} noted that
\begin{quotation}
\noindent
Historically reflection principles are associated with attempts to
formulate the idea that no one notion, idea, statement can capture
our whole view of $V$.
\end{quotation}

It may be surprising to see that we're discussing how to climb higher
in Figure~\ref{fig:largechart} by contemplating ``skipping to the end''
of all possible infinities.
It also sounds somewhat theological, like Anselm's ontological
argument for the existence of God, or Aquinas' argument about
degrees of perfection \citep[see, e.g.,][]{WH16}.
In addition, it may not be surprising to learn that the Reflection
Principle seems to skirt up against the boundary of contradiction!
After all, one ``property'' of $V$ is that it contains all sets.
However, no subset of $V$ (that is not equal to $V$ itself) has
{\em that} property.
The philosopher Bertrand Russell began to study this apparent
paradox at the beginning of the 20th century.
Since then, mathematicians have proven themselves clever enough to
formulate more precise versions of the Reflection Principle that
sidestep the contradiction by specifying their terms in well-formed ways
\citep[see, e.g.,][]{Bg23}.

In any case, the Reflection Principle has been used to help specify
the precise set-theory definitions of just about all of the large
cardinals, from the inaccessibles on up.
However, once we reach the uppermost heights of
Figure~\ref{fig:largechart} and start encountering the
{\em extendible cardinals} (introduced by
\citeauthor{Mag71} [\citeyear{Mag71}] and
\citeauthor{Rein74} [\citeyear{Rein74}]),
I can't help but wonder if {\em some}
aspect of the Reflection Principle starts to break down.
I can only report on the descriptions made by others to state that
an extendible cardinal is so large that it becomes difficult to
distinguish where we even are in the hierarchies of cardinality
and consistency strength.

\citet{Ru82} provided the poetic analogy of climbing up a
steep cliff-face.
As we climb, we can look down to see the lowest levels of cardinality
(i.e., the inaccessible and Mahlo cardinals), and we can see that
there is still a lot of cliff above us, yet to climb.
What would happen if we closed our eyes, and an eagle grabbed us and
flew us up, up, up, and dropped us off at some higher place on the cliff?
If we had started at some middling level (say, at the measurable cardinals),
when we opened our eyes we {\em would} be able to notice the change.
It would be possible to see how much further above the lowest
levels we have flown, and how much closer we are now to visible
landmarks above us.
However, if we had started in the vicinity of the extendible cardinals,
there would be {\em no way} to tell the difference.
The extendible cardinals are so large that virtually everything
surrounding them takes up enough real-estate to make all other
distinctions of size seem insignificant in comparison.

\section{Conclusions}
\label{sec:conc}

Unfortunately, this is where I have to stop.
I really haven't begun to understand many of the large cardinals,
not to mention the huge, superhuge, ultrahuge, and rank-into-rank
cardinals near the top of Figure~\ref{fig:largechart}.
At the {\em very} top sits a statement of utter contradiction,
{\em $0=1$.}
This essentially means that the ``strongest'' possible way to proceed
is to start with an axiom like $0=1$, which would then
mean that we could prove {\em anything we want.}
However, I probably need to paraphrase {\em The Incredibles} and
note that when every possible statement is mutually consistent,
then nothing really is.

I began this journey with appreciation for \citet{Ru82}, and I should
mention that he spent a lot of time discussing how contemplating the
infinite can be considered a union of two opposite modes of thinking.
On one hand, it's the pinnacle of pure MYSTICISM.
In what other mode of life would one dare to approach ``the all?''
We mentioned Cantor's religious reverence of the Absolute $\Omega$,
and I think it bears some resemblance to the concept of 
{\em nirguna brahman} (unmanifest infinite reality) from the
Hindu Upanishads, too.
In addition, recall all the times that we were required to think in
new and creative ways, doing an end-run around some established idea.
This reminds me of how wizards and shamans were said to battle one
another---by reshaping reality to outwit the other---in myths and stories.

On the the other hand, because every step of this effort has been
buttressed by the need to prove things from first principles,
it also seems like the pinnacle of pure RATIONALITY.
This is likely what Hilbert meant by
``the paradise that Cantor has created for us.''
It is an amazing thing that our finite brains have been able to
derive and discern truths about concepts that may in fact be too
big to fit into our physical universe.
The process of learning more and more logical detail about how the
world works should make it all more wondrous,
not less.\footnote{And, by ``the world,'' I mean much more than the
physical universe of atoms and energy.}
That thought may be a clich\'{e}, but others have explored it in
fascinating detail \citep[e.g.,][]{Ei60,Na08,De11}.

Lastly, I need to acknowledge that I haven't really covered
every interesting development that led to the full blooming of the
tree shown in Figure~\ref{fig:largechart}.
I also ended up not mentioning many related topics---e.g., the Cantor
set, the Banach-Tarski paradox, or the intricacies of ordinal
and cardinal arithmetic---that didn't seem to be crucial for
the process of climbing that tree.
My ultimate goal is still to understand and explain all of the
different varieties of infinity for a non-specialist audience,
and I will amend these notes if I find interesting new ways to do that.

\vspace*{0.10in}
\noindent
\hrulefill

\appendix
\section{Proofs That the Rationals Have ``Measure Zero''}
\label{sec:appA}

It turns out that there are several interesting and complementary ways
to show that the {\em countable} set of all rational numbers essentially
takes up zero space along the {\em uncountable} real number line.
This is discussed frequently in textbooks on real analysis
\citep[e.g.,][]{Ru76,Boas96,Ca00} and measure theory \citep{Tao11},
but I have not been able to find a concise collection of these
multiple proofs written for non-specialist audiences.
The subsections below present a few interesting examples.

\subsection{Repeating Versus Non-Repeating Decimals}

This first one is only a plausibility argument, not really a proof.
However, it has the benefit of involving relatively familiar concepts.
Consider the fact that all rational numbers have decimal expansions
that end in an infinitely repeating pattern of digits.
Whether it is only a single repeating digit, such as
\begin{displaymath}
  1 \, = \, 1.000000 \ldots \, = \, 0.999999 \ldots
  \,\,\,\,\,\,\,\,\,\,\,\,\,\,
  \mbox{or}
  \,\,\,\,\,\,\,\,\,\,\,\,\,\,
  \frac{8}{9} \, = \, 0.8888888 \ldots \, = \, 0.\overline{8}
\end{displaymath}
or repeating clumps of multiple digits, like
\begin{displaymath}
  \frac{1}{13} \, = \, 0.076923076923076923076923076923076923 \ldots
  \, = \, 0.\overline{076923} \,\,\, ,
\end{displaymath}
if a number can be represented as the ratio of two natural numbers, then
its decimal expansion eventually terminates in a repeating pattern.

Once we know this, we see that it is possible to {\em deviate} from
this pattern and obtain an irrational number.
We must be careful about how this is done.
If we start with a given rational number and change only one digit,
with all succeeding digits still obeying the original pattern, the
new number is still a rational number.
It is now a different one from the original number, but it's still rational.
However, there are an infinite number of ways of producing
infinitely perturbed (but not regularly-spaced) deviations!
For example, consider something called Liouville's constant, which
is determined by starting with zero and adding 1-digits at non-periodic
intervals determined by the factorial function:
\begin{displaymath}
  x \,\, = \,\, \sum_{n=1}^{\infty} \, \frac{1}{10^{n!}} \,\,\, .
\end{displaymath}
Thus, $x = 0.1 + 0.01 + 10^{-6} + 10^{-24} + 10^{-120} + \cdots$.
Every digit is well-determined, but it does not ever terminate in the
kind of infinitely repeating pattern that defines a rational number.

Now consider generalizing Liouville's constant to something like
\begin{displaymath}
  x(a,b,c) \,\, = \,\, \sum_{n=a}^{\infty} \, \frac{b}{c^{n!}} \,\,\, ,
\end{displaymath}
where $a$, $b$ and $c$ can be any positive integer.
Then, consider adding every possible value of $x(a,b,c)$ to a given
rational number.
This gives us at least $\omega \cdot 3$ different ways of ``perturbing''
that rational number's digits so that it becomes an irrational number.

Of course, there is an even more vastly infinite number of {\em other}
ways of defining a non-periodically perturbed decimal expansion, besides
using the factorial function.
As one more example, consider Champernowne's constant:
\begin{displaymath}
  0.12345678910111213141516171819202122232425 \ldots
\end{displaymath}
i.e., the concatenation of each succeeding decimal numeral in the set
of all natural numbers.
This is irrational, too!
One could also define $\omega$ different unique versions of this kind
of constant, simply by starting at a different natural number (besides 1)
after the decimal point.
Each rational number can be transformed into an infinite number of
new irrationals by adding these constants to it.

Thus, for every rational number, there exist infinitely many irrational
numbers.  If one stakes out some piece of the real line (say, the
interval [0,1] that we examined before), the rational numbers inside
that space must occupy a vanishingly small part of it.

\subsection{Reals and Rationals as Continued Fractions}

Our familiarity with decimal expansions (as the default mode of
expressing real and rational numbers) may just be a historical fluke.
There are other ways of representing numbers to arbitrary precision,
and for this application there is one that may be superior to the
decimal expansion.
Consider the {\em continued fraction,}
\begin{displaymath}
  a_0 +
  \cfrac{1}{a_1 +
  \cfrac{1}{a_2 +
  \cfrac{1}{a_3 +
  \cfrac{1}{ \ddots + \cfrac{1}{a_n} }}}}
\end{displaymath}
which is an iterative way of expressing a number as the sum of
reciprocals that involve only integers $a_n$.
The collection of partial quotients is often expressed using a more
compact notation
\begin{displaymath}
  [ a_0 ; \, a_1, a_2, a_3, a_4, \ldots , a_n ] \,\,\, ,
\end{displaymath}
and each unique collection of $a_n$ quotients has a
one-to-one-correspondence with a specific real or rational number.

The key fact to know is that rational numbers can always be represented
by a {\em finite} number of partial quotients,
but irrational numbers always require an {\em infinite} number
of partial quotients.
Some irrationals (like $\pi$) have no discernible pattern in their
infinite list of partial quotients.
Others (like $\sqrt{2}$) have infinitely repeating digits in their
pattern, and still others (like $e$) have non-periodic patterns
reminiscent of the
decimal representations of Liouville's or Champernowne's constants.
The irrational number with the slowest possible convergence is
the golden ratio,
\begin{displaymath}
  \varphi \,\, = \,\, \frac{1 + \sqrt{5}}{2}  \,\, = \,\,
  [ 1 ; \, 1, 1, 1, 1, 1, 1, \ldots ] \,\,\, .
\end{displaymath}
Thus, there is a clear distinction between a countable
collection of finite sets (which correspond to the
rational numbers), and an uncountably {\em larger} collection of
infinite sets (which correspond to irrational numbers).
Thus, our concept that the rationals occupy much less space than
the irrationals is essentially another way of stating that
$\aleph_0 < {\mathfrak c}$.

\subsection{The Standard Textbook Proof}

This is probably the most frequently seen proof.
First, let us choose a finite piece of the real number line.
The familiar interval [0,1] is a fine choice.
Then, let us create an enumerated list of all rational numbers that
fit into this interval.
This is a similar procedure to the one outlined in
Section~\ref{sec:countable:lim2},
but here we should be careful to exclude any fractions that are divided
by zero, and also exclude duplicates (i.e., once we include 1/2 in the
list, we should skip 2/4, 3/6, 4/8, and so on).
Now that we have an ordered list of rationals, let us call them
$Q_n$, where the subscript $n$ covers all of the natural numbers.

The question now becomes: How much ``real estate'' do these numbers
occupy along the number line?
That is a difficult question to answer because each number is really
just a dimensionless point with no width.
However, let us err on the side of potentially {\em overestimating}
this result by defining some finite widths around each rational number.
Inspired by the speed-up's defined by \citet{Ru82}, let us specify some
finite quantity $\epsilon$ and posit that each rational number $Q_n$
is surrounded by an interval
\begin{displaymath}
  \left[ Q_n - \frac{\epsilon}{2^{n+2}} \,\, , \,\,\,
         Q_n + \frac{\epsilon}{2^{n+2}} \right].
\end{displaymath}
In other words,
$Q_0$ corresponds to a finite line segment of total length $\epsilon / 2$,
$Q_1$ corresponds to a finite line segment of total length $\epsilon / 4$,
$Q_2$ corresponds to a finite line segment of total length $\epsilon / 8$,
and so on.
Given the infinite series of Equation~(\ref{eq:speedup}), we see that
the total length occupied by the full set of rationals is equal to
exactly $\epsilon$.

{\em Aside:} We should admit that the total length may be a bit
less than $\epsilon$ because
the finite interval for any one rational number may overlap a bit with
the finite intervals for some of its neighbors.
But the total length cannot be {\em more} than $\epsilon$.

Note that we never specified a value for $\epsilon$.
All we really know is that it must be a positive real number, so
that each interval has a non-zero length.
The key insight is that the above argument works perfectly well if
we assume that $\epsilon = 0.0000001$, or even if 
$\epsilon = 10^{-999}$.
Thus, it is safe to assume that $\epsilon$ is
infinitesimally small compared to the actual length of the
real interval between 0 and 1.
Thus, because our entire $\epsilon$ construction was always
an overestimate,
it is perfectly fine to conclude that the full set of rational
numbers (of ordinality $\omega^2$, or cardinality $\aleph_0$)
takes up virtually no space along the real number line.
The technical term, which is useful if you want to search for
additional information on this, is that we have shown the rationals
have ``Lebesgue measure zero.''

\end{document}